\documentstyle[epsf]{article}

\newtheorem{theo+}           {Theorem}      [section]
\newtheorem{prop+}  [theo+]  {Proposition}
\newtheorem{coro+}  [theo+]  {Corollary}
\newtheorem{lemm+}  [theo+]  {Lemma}
\newtheorem{exam+}  [theo+]  {Example}
\newtheorem{rema+}  [theo+]  {Remark}
\newtheorem{defi+}  [theo+]  {Definition}

\newtheorem{exam+s}  [theo+]  {Examples}
\newtheorem{rema+s}  [theo+]  {Remarks}
\newtheorem{hyp+}  [theo+]  {Hypotheses}
\newtheorem{cla+}  [theo+]  {Claim}

\newenvironment{theorem}{\begin{theo+}}{\end{theo+}}
\newenvironment{proposition}{\begin{prop+}}{\end{prop+}}
\newenvironment{corollary}{\begin{coro+}}{\end{coro+}}
\newenvironment{lemma}{\begin{lemm+}}{\end{lemm+}}
\newenvironment{example}{\begin{exam+}\rm}{\end{exam+}}
\newenvironment{remark}{\begin{rema+}\rm}{\end{rema+}}

\newenvironment{examples}{\begin{exam+s}\rm}{\end{exam+s}}
\newenvironment{remarks}{\begin{rema+s}\rm}{\end{rema+s}}

\newcommand{\pa}{\partial}
\newcommand{\half}{{\textstyle\frac{1}{2}}}
\newcommand{\be}{\begin{equation}}
\newcommand{\ee}{\end{equation}}
 
\newcommand{\Bbb}{\bf}
 
\newcommand{\RR}{{\Bbb R}}  
\newcommand{\CC}{{\Bbb C}}  
\newcommand{\HH}{{\Bbb H}}  
\newcommand{\MM}{{\Bbb M}}  

\newcommand{\Ff}{{\cal F}}  
\newcommand{\Ll}{{\cal L}}  
\newcommand{\Cc}{{\cal C}}  
\newcommand{\Jj}{{\cal J}}  

\newcommand{\Rt}{\widetilde{\Bbb R}^3}  
\newcommand{\Rf}{\widetilde{\Bbb R}^4}  
\newcommand{\Cf}{\widetilde{\Bbb C}^4}  
\newcommand{\Mf}{\widetilde{\Bbb M}^4}  
\newcommand{\Ci}{{\Bbb C} \cup \{ \infty\} }  

\newcommand{\CP}[1]{{{\Bbb C}P}^{#1}} 

\newcommand{\ii}{{\rm i}}  
\newcommand{\jj}{{\rm j}}  

\newcommand{\ra}{\rightarrow}
\newcommand{\lra}{\longrightarrow}
\newcommand{\llra}{\longleftrightarrow}
\newcommand{\proof}{\noindent\ {\bf Proof} \hskip 0.4em}
\newcommand{\eproof}{\bigskip}

\begin{document} 
\title{Harmonic morphisms, conformal foliations and shear-free ray
congruences}

\author{P. Baird and J. C. Wood
\thanks{Partially supported by EC grant CHRX-CT92-0050} \\ {\small
D\'epartement de Math\'ematiques, Universit\'e de Bretagne Occidentale} \\
{\small 6 Avenue Le Gorgeu, B.P. 452, 29275 Brest Cedex, France, and} \\
{\small Department of Pure Mathematics, University of Leeds}\\ {\small
Leeds LS2 9JT, G.B.} }

\date{} 
\maketitle

\begin{abstract}
Equivalences between conformal foliations on Euclidean $3$-space,
Hermitian structures on Euclidean $4$-space, shear-free ray
congruences on Minkowski $4$-space, and holomorphic foliations on
complex $4$-space are explained geometrically and twistorially;
these are used to show that 1) any real-analytic complex-valued
harmonic morphism without critical points defined on an open
subset of Minkowski space is conformally equivalent to the
direction vector field of a shear-free ray congruence, 2) the
boundary values at infinity of a complex-valued harmonic morphism
on hyperbolic $4$-space define a real-analytic conformal
foliation by curves of an open subset of Euclidean $3$-space and
all such foliations arise this way.  This gives an
explicit method of finding such foliations; some examples are given.
\end{abstract}

\section{Introduction}
This paper consists of two parts:  Firstly, we describe natural
correspondences between the following four quantities:

(Q1) holomorphic foliations by $\alpha$-planes of an open subset
of complex $4$-space $\CC^4$,

(Q2) positive Hermitian structures $J$ on an open subset of
Euclidean $4$-space $\RR^4$,

(Q3) shear-free ray congruences on an open subset of Minkowski
$4$-space $\MM^4$,

(Q4) conformal foliations by curves of an open subset of
Euclidean $3$-space $\RR^3$.

The correspondences are described both geometrically and
twistorially.

Then these correspondences are used to find a relationship between
(complex-valued) {\em Minkowski\/} harmonic morphisms, i.e. ones
defined on open subsets of Minkowski space, and shear-free ray
congruences; in fact, the direction vector field of a shear-free
ray congruence defines such a harmonic morphism; conversely, to
any Minkowski harmonic morphism $\phi$ without critical points,
we can associate a shear-free ray congruence such that every
fibre of $\phi$ is the union of parallel rays of the congruence
(Theorem \ref{th:SFR}).  It follows that, up to conformal
transformations of the codomain, any ``non-K\"ahler"
Minkowski harmonic morphism
is the direction vector field of a shear-free ray congruence
(Corollary \ref{cor:compo}). Similar results are given for
{\em complex\/} harmonic morphisms (Theorem \ref{th:cxhamorphSFR}
and Corollary \ref{cor:compo}).

Secondly we show (Theorem \ref{th:HWC-hyp}) that the level sets
of the boundary values at infinity of a complex-valued submersive
harmonic morphism from hyperbolic $4$-space define
a conformal foliation of an open subset of $\RR^3$ by curves and that
every  such (real-analytic) foliation arises in this
way. Combining with \cite{Wood-4d} and \cite{Baird-JMP}, this
gives us a practical method of finding all conformal foliations
by curves of open subsets of $\RR^3$; this is illustrated in the
final chapter where some explicit examples are calculated.

We now describe these ideas in more detail:

A {\em harmonic morphism\/} is a smooth map $\phi$ between Riemannian
manifolds which preserves Laplace's equation in the sense that if
$f$ is a local harmonic function on the codomain manifold then $f
\circ \phi$ is a local harmonic function on the domain, see
\cite{Baird-JMP} and \cite{Wood-4d} for some relevant background.

{\em A conformal foliation by curves\/} (\S \ref{subsec:conf}) of
an open subset of Euclidean $3$-space $\RR^3$ is a foliation by
curves such that the Lie transport along the fibres is conformal
on the normal bundle.  Conformal foliations include Riemannian
foliations.  A foliation is conformal if and only if its leaves
are locally the level sets of a complex-valued submersion $f$
which satisfies the {\em horizontal weak conformality
condition\/}:
\be
\left( \frac{\pa f}{\pa x_1} \right)^2 + \left(\frac{\pa f}{\pa
x_2} \right)^2 + \left( \frac{\pa f}{\pa x_3 } \right)^2 = 0\;.
\label{HC0} \ee

{\em A positive Hermitian structure\/} $J$ (\S
\ref{subsec:distHerm}) on an open subset of Euclidean $4$-space
$\RR^4$ is an almost Hermitian structure which is associated to
a positive basis and is integrable.

{\em A shear-free ray congruence} (\S \ref{subsec:SFR})
on an open subset of
Minkowski $4$-space $\MM^4$ is a foliation by null lines such
that Lie transport is conformal on the {\em screen spaces\/}.

{\em A holomorphic foliation by $\alpha$-planes} (\S
\ref{subsec:distnull}) of an open subset of complex $4$-space
$\CC^4$ is a foliation, holomorphic with respect to the standard
complex structure on $\CC^4$, by planes which are null and are
associated to a positively oriented basis.

We show that these quantities are equivalent by establishing
correspondences between them.  For
example:

Given a Hermitian structure $J$ on an open subset of $\RR^4 =
\{(x_0,x_1,x_2,x_3)\}$, putting $U = J(\pa/\pa x_0)$ on a slice
$x_0 =$ constant defines the tangent vector field of a
real-analytic conformal foliation by curves of an open subset of
$\RR^3$ and any such foliation arises this way giving a
bijective correspondence between germs (Corollary \ref{cor:JtoU}).

Analogously, given a $(C^{\infty})$ shear-free ray congruence on
an open subset of $\MM^4  = \{(t,x_1,x_2,x_3)\}$, projecting its
tangent vector onto a slice $t =$ constant defines the tangent
vector field of a $(C^{\infty})$ conformal foliation by curves of
an open subset of $\RR^3$; any such foliation arises this way
giving a bijective correspondence between germs (Theorem
\ref{th:SFR-conf}). In fact, as $t$ varies, we get a
$1$-parameter family of ``associated" conformal foliations.

Both shear-free ray congruences and Hermitian structures
complexify to holomorphic foliations by $\alpha$-planes (see
Proposition \ref{prop:unifthree}) so that all four quantities are
in bijective correspondence (Theorem \ref{th:unifgerms}) and have
a Kerr-type representation as a complex hypersurface of the
twistor space $\CP{3}$ (\S \ref{subsec:Kerr}).

The quantities (Q1 - Q4) arise as distributions of the quantities
at a point in (\ref{QP}) below which satisfy certain conditions;
we first show that these quantities are equivalent at a point
(Proposition \ref{prop:bijpt}) and can all be represented by a
unit vector in $\RR^3$, {\em the direction vector\/} of the
quantity, or twistorially (\S 2.7), this last description
extending to compactified spaces. Then in \S 3 we consider how
various conditions on distributions of the quantities (\ref{QP})
correspond leading to the equivalences of (Q1 - Q4). We spend
some time explaining this material as, although some aspects may
be known to Mathematical Physicists, we hope that it is
worthwhile to give a geometrical (index and spinor-free)
description adapted to our applications. In \S 5 we give our main
result (Theorem \ref{th:HWC-hyp}) and its application to finding
conformal foliations of open subsets of $\RR^3$ by curves. The
idea is that, given a complex hypersurface $S$ of $\CP{3}$, $S$
is the twistor surface of some holomorphic structure $J$ and we
can find a hyperbolic harmonic morphism $\phi$ holomorphic with
respect to $J$ by solving a first order holomorphic partial
differential equation (\ref{E}) on $S$.  Then the level sets of
the boundary values on the $\RR^3$ at infinity define a
conformal foliation.  By introducing a parameter $a \in \CC^4$ we can
actually find a $5$-parameter family of associated conformal
foliations.

Some of our work can be generalized to more
general manifolds but with a loss of explicitness; this will be
done in a future article.

\bigskip {\bf Acknowledgments}

The authors would like to thank U. Pinkall for a suggestion which
led to this investigation and T. Bailey  for conversations
related to this work. The second author would like to thank the
D\'epartement de Math\'ematiques at the Universit\'e de Bretagne
Occidentale, Brest, for making possible a visit in July 1995
where much of this work was done.


\section{Null planes and associated structures}
In this section we consider four quantities defined at a point $p$ of
$\CC^4$ or its conformal compactification and show that they are
all in bijective correspondence, namely, with definitions to
follow,
\be \left. \begin{array}{rl}
\mbox{(QP1)} & \alpha \mbox{-planes } \Pi_p \mbox{ at } p\,, \\
\mbox{(QP2)} & \mbox{positive real almost Hermitian structures }
J_p \mbox{ at } p\,, \\
\mbox{(QP3)} & \mbox{null directions } V_p \mbox{ at } p
\mbox{ in the Minkowski slice through } p\,, \\
\mbox{(QP4)} & \mbox{unit vectors } U_p \mbox{ at } p
\mbox{ in the } \RR^3
\mbox{-slice through } p\,.
\end{array} \right\} \label{QP} \ee

\subsection{Null planes} \label{Null planes}

We consider $\CC^4 = \{(x_0,x_1,x_3,x_4):x_i \in \CC \}$ with its
standard complex structure, and conformal structure given by the
holomorphic metric $g = \sum_{i=0}^3 dx_i^{{}2}$. (For general
definitions of holomorphic metrics and related concepts, see
\cite{LeBruncx}.)  A tangent vector $V = (V_0,V_1,V_2,V_3)
\in T_p\CC^4,\ p \in \CC^4$ is called {\em null\/} if $g(V,V) =
0$.  A {\em null plane at } $p \in \CC^4$ is a two-dimensional
complex linear subspace of $T_p\CC^4$ consisting of null vectors.
 An {\em affine null plane in } $\CC^4$ is a $2$-dimensional
complex affine subspace of $\CC^4$ whose tangent space at any
point is a null plane.

For any $p \in \CC^4$, let $\RR^4_p$ be the real $4$-dimensional
affine subspace ({\em real slice\/}) through $p = (p_0, p_1, p_2,
p_3)$ given by $\{ (x_0,x_1,x_2,x_3) \in \CC^4 : \Im x_i = \Im
p_i,\ i=0,1,2,3 \}$ and parametrized by $\RR^4 \ni
(x_0,x_1,x_2,x_3) \mapsto (p_0 + x_0, p_1 + x_1, p_2 + x_2, p_3 +
x_3) \in \CC^4$; if $p = 0 \equiv (0,0,0,0)$ we write $\RR^4_0 =
\RR^4$.  We give each $\RR^4_p$ the standard orientation.  The
holomorphic metric and conformal structure restrict to the
standard metric and conformal structure on each $\RR^4_p$.  Note
that $T_p \CC^4$ can be identified with $T_p \RR^4_p \otimes
\CC$; we shall frequently make this identification.

Given any orthonormal basis $\{ e_0,e_1,e_2,e_3 \}$ of $T_p
\RR^4$, the plane $\Pi_p = \mbox{span}\{e_0+\ii e_1, e_2 + \ii
e_3\}$ is null; we call $\Pi_p$ an {\em $\alpha$-plane\/} (resp.
{\em $\beta$-plane\/}) according as the basis $\{ e_0,e_1,e_2,e_3
\}$ is positively (resp. negatively) oriented.  This construction
induces a bijection
\be
\mbox{S0(4)/U(2)} \llra \{ \alpha \mbox{-planes at } p \}\;.
\label{alpha} \ee

To proceed, it is convenient to introduce new coordinates $(z_1,
\widetilde{z}_1, z_2, \widetilde{z}_2)$ on $\CC^4$ by setting 
\be
z_1 = x_0 + \ii x_1, \ \widetilde{z}_1 = x_0 - \ii x_1, \ z_2 = x_2 +
\ii x_3, \ \widetilde{z}_2 = x_2 - \ii x_3.
\label{coords} \ee
Then $\RR^4$ is given by $\widetilde{z}_1 = \overline{z}_1, \widetilde{z}_2 =
\overline{z}_2$ and the holomorphic metric on $\CC^4$ by $g = dz_1
d\widetilde{z}_1 + dz_2 d\widetilde{z}_2$.

\subsection{Real almost Hermitian structures and null-planes}
By a {\em real almost Hermitian structure\/ $J_p$ at\/ $p \in
\CC^4$} we mean an endomorphism $J_p:T_p\CC^4 \ra T_p\CC^4$ which
maps real vectors, i.e. ones in $T_p\RR^4_p\,$, to real vectors
isometrically and has $J_p^2 = -I$.  Equivalently, a real almost
Hermitian structure is an almost Hermitian structure on the real
slice $\RR^4_p$ at $p$ extended to $T_p \CC^4 = T_p \RR^4_p
\otimes \CC$ by complex linearity.

Given any orthonormal basis $\{ e_0, e_1, e_2, e_3 \}$ of $T_p
\RR^4_p$, setting $J_p(e_0) = e_1,\ J_p(e_2) = e_3$ defines a
real almost Hermitian structure $J_p$ at $p$; we call $J_p$ {\em
positive\/} (resp. {\em negative\/}) according as $\{
e_0,e_1,e_2,e_3 \}$ is a positively (resp. negatively) oriented
basis; this construction gives a bijection
\be
\mbox{SO(4)/U(2)} \llra \{ \mbox{positive real almost Hermitian
structures } J_p \mbox{ at } p \}\;. \label{H}
\ee
Combining with (\ref{alpha}) gives a bijection
$$
\{ \alpha \mbox{-planes } \Pi_p \mbox{ at } p \} \llra \{
\mbox{positive real almost Hermitian structures } J_p \mbox{ at }
p \} $$
given by
$$
\Pi_p \mapsto J_p = \left\{ \begin{array}{l} -\ii \mbox{ on }
\Pi_p \\ +\ii \mbox{ on } \overline{\Pi_p} \end{array} \right.
$$
with inverse
$$
J_p \mapsto \Pi_p = (0,1)\mbox{-tangent space of } J_p = \{ X+\ii
J_p X : X \in T_p \RR^4_p \}.
$$

\subsection{Vectors in $\RR^3$}
For any $p = (p_0,p_1,p_2,p_3) \in \CC^4$, define the {\em
$\RR^3$-slice through $p$} by
$$
\RR^3_p = \{ (x_0,x_1,x_2,x_3): x_0 = p_0, \Im x_i = \Im p_i \
(i=1,2,3) \}
$$
parametrized by $\RR^3 \ni (x_1,x_2,x_3) \mapsto (p_0,p_1 + x_1,
p_2 + x_2, p_3 + x_3) \in \CC^4$. Given any unit vector $U_p \in
T_p \RR^3_p$, there is a positive orthonormal basis $\{e_0, e_1,
e_2, e_3 \}$ of $T_p\RR^4_p$ with $e_0 = \pa / \pa x_0 ,\ e_1 =
U_p\,$, so that $\{ e_2, e_3 \}$ spans $U_p^{\perp} \cap
T_p\RR^3_p\,$.  This then defines an $\alpha$-plane $\Pi_p =
\mbox{span} \{e_0+\ii e_1, e_2+ \ii e_3 \}$ and so a positive
real almost Hermitian structure $J_p$, and gives rise to a
bijection:
\be
\{ \mbox{positive real almost Hermitian structures } J_p \mbox{ at
} p \} \leftrightarrow
\{ \mbox{unit vectors } U_p \mbox{ in } T_p\RR^3_p \} 
\label{JU}\ee
given by: $J_p \mapsto U_p = J_p(\pa / \pa x_0)$, \ $U_p \mapsto
J_p =$ the unique positive almost Hermitian structure at $p$ with
$J_p(\pa / \pa x_0) = U_p$.

\subsection{Null vectors in $\MM^4$ and the fundamental
bijections at a point}
Minkowski space $\MM^4$ is defined to be the set $\{ (t, x_1, x_2,
x_3) \in \RR^4 \}$ with the semi-Riemannian metric $g^M = -dt^2 +
dx_1^{{}2} + dx_2^{{}2} + dx_3^{{}2}$.  We include $\MM^4$ in
$\CC^4$ by the map $(t, x_1, x_2, x_3) \mapsto (-\ii t, x_1, x_2,
x_3)$.  (Note that the minus sign is unimportant; it is included
simply to avoid minus signs later on.)  More generally, for any
$p \in \CC^4$, let $\MM^4_p$ be the {\em Minkowski slice\/}
$\MM^4_p = \{ (x_0, x_1, x_2, x_3 ) : \Re x_0 = \Re p_0, \Im x_i
= \Im p_i \ (i=1,2,3) \}$ parametrized by $\MM^4 \ni
(t,x_0,x_2,x_2) \mapsto (p_0 - \ii t, p_1 + x_1, p_2 + x_2, p_3 +
x_3) \in \CC^4$.   A vector $v = (v_0, v_1, v_2, v_3) \in T_p
\MM^4_p$ is called {\em null\/} if $g_M(v,v) = -v_0^{{}2} +
v_1^{{}2} + v_2^{{}2} + v_3^{{}2} = 0$, equivalently its image
$(-\ii v_0, v_1, v_2, v_3)$ in $T_p \CC^4$ is null in the sense
of \S \ref{Null planes}.  A $1$-dimensional subspace $V_p$ of
$T_p \MM^4_p$ is called a {\em null direction  at\/} $p$ if it is
spanned by a null vector $v_p$.  An {\em (affine) null line\/} or
{\em(light) ray\/} of $\MM^4_p$ is a line in $\MM^4_p$ whose
tangent space at any point is null, i.e. spanned by a null
vector.  Any null plane $\Pi_p \subset T_p\CC^4_p$ at $p$
intersects $T_p \MM^4_p$ in a null direction $V_p$ at $p$.
Indeed, if we write $\Pi_p = \mbox{span} \{ e_0 + \ii e_1 , e_2 +
\ii e_3 \}$ where $\{e_0, e_1, e_2, e_3 \}$ is a positive
orthonormal basis of $T_p\RR^4_p$ with $e_0 = \pa/ \pa x_0$, then
$V_p$ is spanned by $v_p = \pa / \pa t + e_1$.  This gives a
bijection between null directions at $p$ and $\alpha$-planes at
$p$ which, combined with our previous bijections gives
correspondences which will be fundamental in this paper:

\begin{proposition}\label{prop:bijpt}
For any $p \in \CC^4$ we have bijective correspondences between
the quantities (QP1) --(QP4) in
 (\ref{QP}) given by

\be \left. \begin{array}{rclcl}
J_p  & = & J_p(\Pi_p) & = & \left\{ \begin{array}{l} -\ii
\mbox{ on } \Pi_p \\
 +\ii \mbox{ on } \overline{\Pi_p} \end{array} \right.
\\
\Pi_p & = & \Pi_p(J_p) & = & (0,1)\mbox{-tangent space of } J_p =
\{ X+\ii J_p X : X \in T_p \RR^4_p \} \\             
U_p  & = & U_p(J_p) & = & J_p(\pa / \pa x_0)  \\
J_p & = & J_p(U_p) & = & \mbox{the unique positive almost Hermitian
structure at }p \\
& & & &\mbox{ with } J_p(\pa / \pa x_0) = U_p  \\ 
V_p  & = & V_p(U_p) & = & \mbox{span}(\pa / \pa t + U_p) \\
U_p & = & U_p(V_p) & = & \mbox{normalized projection of } V_p
\mbox{ onto } T_p\RR^3_p \\
& & & & \mbox{ i.e. the unique vector }U_p \mbox{ such that }V_p =
\mbox{span}(\pa / \pa t + U_p)\;.  
\end{array} \right\} \label{Ufor} \ee
\end{proposition}

Schematically we have a commuative diagram:

\small\be
\begin{array}{ccccc} 
 & & \vbox{ \hbox{\space \space \space \space \space \space \space 
\space \space \space \space \space \space (QP2) $=$}
\hbox{ \{positive real almost Hermitian} \hbox{ \space
\space \space \space \space \space \space \space \space \space
structures at $p$:}
\hbox{\space \space \space \space \space \space \space
$J_p:T_p\RR^4_p \ra T_p\RR^4_p$
\} } } & &  \\ 
 & \nearrow &   & \searrow &   \\
\vbox{ \hbox{\space \space \space \space \space (QP1) $=$}
\hbox{ \{$\alpha$-planes at $p$:} \hbox{\space \space
$\Pi_p \subset T_p\CC^4$
\} } } & & & & \vbox{ \hbox{ \space \space \space (QP4) $=$}
\hbox{  \{unit vectors } \hbox{\space \space
$U_p$ in $T_p\RR^3_p$
\} } }  \\
 & \searrow &   & \nearrow &   \\
 &          &  \vbox{ \hbox{\space \space \space \space \space \space
\space \space (QP3) $=$}
\hbox{ \{null directions at $p$:} \hbox{\space
\space \space \space \space \space 
$V_p \subset T_p\MM^4_p$
 \} } }      &          &
\end{array}
\label{Upt} \ee \normalsize

\noindent
Given any one of the four quantities (\ref{Upt}) at a point $p$ of
$\CC^4$, the relations (\ref{Ufor}) determine the other three
corresonding quantities.  In particular, all four quantities are
determined by a unit vector $U_p \in T_p^1\RR^3_p \cong S^2$, which
we shall call the {\em direction vector\/} of the quantity.

\begin{remark}
There is a fifth quantity in bijective correspondence with the
quantities (\ref{QP}) namely {\em null directions in
$T_p\CC^3_p$} where $\CC^3_p$ denotes the $\CC^3$-slice
$\{(x_0,x_1,x_2,x_3) \in \CC^4 : x_0 = p_0 \}$.
Indeed we have a bijection:
$$
\{ \mbox{unit vectors in } T_p\RR^3_p \} \llra \{ \mbox{null
directions in } T_p\CC^3_p \}
$$
given, in the notations above, by
$$
U_p \llra \mbox{span}\{e_2 +\ii e_3 \} \,.
$$
However, we shall not make use of this quantity.
\end{remark} 

\subsection{Representations in coordinates} \label{subsec:repcoord}
For any $p \in \CC^4$ and $[w_0, w_1] \in \CP{1}$, set
\be
\Pi_p = \mbox{span} \left\{w_0 \frac{\pa}{\pa \widetilde{z}_1}
- w_1 \frac{\pa}{\pa z_2},
w_0 \frac{\pa}{\pa \widetilde{z}_2} + w_1 \frac{\pa}{\pa z_1} \right\}
\;. \label{NT0} \ee
Writing $\mu = w_1/w_0 \in \Ci$, if $\mu \neq \infty$, \  $\Pi_p$
has a basis
\be
\left\{ \frac{\pa}{\pa \widetilde{z}_1} - \mu \frac{\pa}{\pa z_2},
\qquad
 \frac{\pa}{\pa \widetilde{z}_2} + \mu \frac{\pa}{\pa z_1} \right\}
\;; \label{NT}
\ee
if $\mu = \infty$, a basis for $\Pi_p$ is given by
$\left\{ \pa/\pa z_2 , \pa/\pa z_1 \right\}$.

$\Pi_p$ is easily seen to be an $\alpha$-plane at $p$ and
(\ref{NT0}) gives an explicit bijection
$$
\CP{1} \llra \{ \alpha\mbox{-planes at } p \}.
$$
Applying the matrix $\left(
\begin{array}{ll} \overline{w}_0 & \overline{w}_1 \\
                           -w_1  &           w_0  
\end{array} \right)$,
a multiple of which is unitary, converts the basis of $\Pi_p$ to
an equivalent basis $\{e_0+\ii e_1, e_2 + \ii e_3 \}$. Writing $u
= iw_1/w_0 = i\mu$, and letting $\sigma:S^2 \ra \Ci$ denote
stereographic projection from $(-1,0,0)$, the positive
orthonormal basis $\{ e_0, e_1, e_2, e_3 \}$ is given by

\be
\left. \begin{array}{rcl}
e_0 & = & \pa / \pa x_0 \\
U_p = e_1 & = & \displaystyle{\frac{1}{1+|u|^2}} (1-|u|^2, 2\Re u,
2\Im u) = \sigma^{-1}(u) \\
e_2 + \ii e_3 & = & \displaystyle{\frac{1}{1+|u|^2}} \big(-2u,
1-u^2, \ii (1+u^2) \big) \;.
\end{array} \right\}
\label{stdbasis} \ee

Note that $\{ e_2,e_3 \}$ gives an orthonormal basis for the {\em
screen space\/} $U_p^{{}\perp} = V_p^{{}\perp} \cap \RR^3_p$ of
the corresponding null direction \be
V_p = \mbox{span}(\frac{\pa}{\pa t} + U_p)
\label{VPi} \ee
and the corresponding almost Hermitian structure $J_p$ is
determined by
\be
J_p (e_0) = e_1, \quad J_p (e_2) = e_3 \;.
\label{Herm}\ee

In summary, {\em any of the corresponding quantities $\Pi_p \in$
(QP1), $J_p \in$ (QP2), $V_p \in$ (QP3), $U_p \in$ (QP4) can be
represented by a point $[w_0,w_1] \in \CP{1}$ or a number $\mu =
w_1/w_0 \in \Ci$, the quantities being given by (\ref{NT0}) --
(\ref{Herm})}.

\subsection{Compactifications} \label{subsec:compact}
We form the conformal compactification $\Cf$ of $\CC^4$ as follows
(cf. \cite{Ward-Wells}  but note that our conventions differ
slightly from theirs --- they agree with those of
\cite{Wood-4d}):  Set $\Cf = G_2 (\CC^4)$, the Grassmannian of
all complex $2$-planes through the origin with its standard
complex structure.  Then $\CC^4$ is embedded holomorphically in
$\Cf$ by the mapping \be
\iota : (z_1, \widetilde{z}_1, z_2, \widetilde{z}_2) \mapsto \mbox{ column
space of } \left[\begin{array}{rr} 1\ & 0\ \\ 0\ & 1\ \\ z_1  &
-\widetilde{z}_2 \\ z_2 & \widetilde{z}_1 \end{array} \right].
\label{emb} \ee
We shall call this the {\em standard coordinate chart\/} for
$\Cf$.  Points of $\Cf \setminus \CC^4$ will be called {\em
points at infinity}. This embedding gives the image of $\CC^4$ in
$\Cf$ a holomorphic conformal structure; this can be extended to
the whole of $\Cf$ by the action of $SU(2,2)$ on $\Cf$;  this
group restricting to the conformal group on $\CC^4$ which sends
(affine) null lines to null lines.

Real slices $\RR^4_p$ compactify to submanifolds $\Rf_p$
conformally equivalent to $4$-spheres; for example the map $\RR^4
= \RR^4_0 \hookrightarrow \Cf$ is given by sending $(z_1,
\overline{z}_1, z_2, \overline{z}_2) \in \RR^4 \mapsto
\mbox{column space of }\displaystyle{\left[ \frac{I}{Q} \right]}$
where $Q =  \left[ \begin{array}{rr}  z_1 & -\overline{z}_2 \\
z_2 & \overline{z}_1 \end{array} \right]$ is the usual
representation of the quaternion $z_1 + z_2 \jj\,$; then $\RR^4$
is compactified by adding the single point at infinity given by
the column space of $\displaystyle{\left[ \frac{0}{I} \right]}$.

Each $\RR^3$-slice $\RR^3_p$ compactifies in $\Cf$ to a
submanifold $\Rt_p$ conformally equivalent to the $3$-sphere
$S^3$.

Each $\MM^4_p$ compactifies in $\Cf$ to a subset $\Mf_p$
topologically equivalent to $S^3 \times S^1$.  For more details,
see, for example, \cite{H-T,P-R,Ward-Wells}.

All the quantities (\ref{QP}) make sense on the compactified
spaces and we still have the bijections (\ref{Upt}), though to
represent them by a vector in $S^2$ requires coordinates to be
chosen.

\subsection{Representations by twistors}
\label{subsec:twistors}
We recall the twistor correspondence which parametrizes the set of
affine null planes in $\Cf$ by $\CP{3}$ (cf. \cite{Ward-Wells}):

Let $F_{1,2}$ be the complex manifold $\{ (w,p) \in \CP{3} \times
\Cf :$ the line $w$ lies in the plane $p \}$, the
so-called {\em correspondence space\/}.  The restrictions of the
natural projections define the double holomorphic
fibration:
\be
\begin{array}{ccccc} 
       &              & (w,p) \in F_{1,2} &              &   \\
       & \mu \swarrow &                   & \searrow \nu &   \\
w \in \CP{3} &        &                &                 & p \in \Cf
\end{array}
\label{doublefibr} \ee

Then for any $w \in \CP{3}$, $\widetilde{w} = \nu \circ \mu^{-1}(w)$
is an $\alpha$-plane in $\Cf$ which we call the $\alpha$-plane
{\em determined by\/} or {\em represented by\/} $w$.

Conversely, for any point $p \in \Cf$ we write $\widehat{p} = \mu
\circ \nu^{-1}(p)$; then $\widehat{p}$ is a $\CP{1}$ in $\CP{3}$
which represents the $\CP{1}$'s-worth of $\alpha$-planes through
$p$.

Explicitly, in the standard coordinates (\ref{emb}) for $\CC^4
\hookrightarrow \Cf$, $(w,p) \in F_{1,2}$ if and only the
following {\em incidence relations\/} are satisfied:
\be
w_0 z_1 - w_1 \widetilde{z}_2 = w_2 \qquad w_0 z_2 + w_1 \widetilde{z}_1 =
w_3
\label{I} \ee
(for these express the condition that $w$ is a linear combination
of the columns of the matrix (\ref{emb}) and so lies in the plane
represented by the point $p \in G_2(\CC^4)\;$).  Now note that,
for any $[w_0,w_1,w_2,w_3]$ with $[w_0,w_1] \neq [0,0]$,
(\ref{I}) defines an affine $\alpha$-plane in $\CC^4$: indeed,
its tangent space at any point is the set of
vectors annihilated by $\{w_0 dz_1 - w_1 d\widetilde{z}_2, w_0 dz_2 +
w_1 d\widetilde{z}_1 \}$ and so is given by (\ref{NT0}).  Points of
$\CP{3}$ on the projective line $\CC P^1_0 = \{ [w_0,w_1,w_2,w_3]
:[w_0, w_1] = [0,0] \}$ correspond to affine $\alpha$-planes at
infinity, i.e. in $\Cf \setminus \CC^4$.  Thus $\CP{3}$
parametrizes all $\alpha$-planes in $\Cf$ and {\em (\ref{I})
expresses the condition that a point $(z_1, \widetilde{z}_1,z_2,
\widetilde{z}_2) \in \CC^4$ lies on the $\alpha$-plane determined by
$[w_0,w_1,w_2,w_3] \in \CP{3}$}, i.e. (\ref{I}) is the equation
of that $\alpha$-plane.

Note that any $\alpha$-plane intersects a slice $\Rf_p$ in a unique
point giving a map $\pi_p:\CP{3} \ra \Rf_p\,$; for example, if
$[w_0,w_1,w_2,w_3] \in \CP{3} \setminus \CC P^1_0$, the
$\alpha$-plane (\ref{I}) intesects $\RR^4 = \RR^4_0$ when
\be
\left. \begin{array}{lll}
                      w_0 z_1 - w_1 \overline{z}_2 & = & w_2 \\
\overline{w}_1 z_1 + \overline{w}_0 \overline{z}_2 & = &
\overline{w}_3 \end{array}
\right\}
\label{IR} \ee
which has the unique solution
\be 
\left( \begin{array}{l} z_1 \\ z_2 \end{array} \right) =
\frac{1}{|w_0|^2 + |w_1|^2}
\left( \begin{array}{l} \overline{w}_0 w_2 + w_1 \overline{w}_3 \\
                \overline{w}_0 w_3 - \overline{w}_1 w_2 \end{array}
\right) \;.
\label{int} \ee
The resulting map
\be
\pi = \pi_0:\CP{3} \ra \Rf = S^4 \label{twistormap}
\ee
is just the standard twistor map --- in quaternionic notation
$[w_0,w_1,w_2,w_3] \mapsto (w_0 + w_1 \jj)^{-1}(w_2+w_3 \jj) \in
\HH \cup \infty \cong S^4$.  The affine $\alpha$-plane (\ref{I})
intersects the Minkowski slice $\MM^4_p$ if and only if it
intersects $\RR^3_p$;  the intersection is then an affine null
line. For $p=0$ this holds if and only if the point (\ref{int})
lies in $\RR^3$, \ i.e. $\Re z_1 = 0$, i.e. $[w_0,w_1,w_2,w_3]$
lies on the real quadric
\be
N^5 = \pi^{-1}(\RR^3) =
\{ [w_0,w_1,w_2,w_3] : w_0 \overline{w}_2 + \overline{w}_0
w_2 + w_1 \overline{w}_3 + \overline{w}_1 w_3 = 0 \} \subset
\CP{3} \;.
\label{N5} \ee
Points in $N^5$ thus represent affine null lines of $\MM^4$.
For general $p$ we replace $N^5$ by $N^5_p = \pi_p^{-1}(\RR^3_p)$.

Interchanging $z_2$ and $\widetilde{z}_2$ in (\ref{I}) gives the
standard parametrization of $\beta$-planes by $\CP{3}$.

The incidence relations (\ref{I}) define a fundamental map which
gives the point of $\CP{3}$ representing the $\alpha$-plane
through $(z_1, \widetilde{z}_1, z_2, \widetilde{z}_2)$ with direction
vector $\sigma^{-1}(\ii w_1/w_0)$: \be
\begin{array}{rrl} \CC^4 \times \CP{1} & \lra   & \CP{3} \setminus
{\CC P^1_0} \ \subset \ \CP{3}    \\
  \big( (z_1, \widetilde{z}_1, z_2, \widetilde{z}_2), [w_0,w_1] \big)
& \mapsto & [w_0, w_1, w_0 z_1 - w_1 \widetilde{z}_2, w_0 z_2 + w_1
\widetilde{z}_1] \;. \end{array}
\label{fund} \ee

The restriction of this to $\RR^4 \times \CP{1}$ gives a
trivialisation of the twistor bundle (\ref{twistormap}) over
$\RR^4$.

In summary, {\em (i) $w \mapsto \widetilde{w}$ defines a bijection
from $\CP{3}$ to the set of all affine $\alpha$-planes in $\Cf$;
if $w \in \CP{3} \setminus \CP{1}_0$, \ $\widetilde{w}$ is given on
$\CC^4$ by (\ref{I}), (ii) a point $(w,p) \in F_{1,2}$ represents
the $\alpha$-plane $\Pi_p$ at $p$ tangent to $\widetilde{w}$ and,
so, any of the four quantities (\ref{QP}); at a point $p \in
\CC^4$, $\Pi_p$ is given by (\ref{NT0})}.


\section{Distributions of null planes and associated
structures}
In this section we study how our four quantities (\ref{QP}) vary
when the point is varied and obtain bijections between their
germs (Theorem \ref{th:unifgerms}).  As in the last section, for
an open set $A^C$ of $\CC^4$,\, $g$ will denote the standard
holomorphic metric or any conformal multiple of it, and $\nabla$
its Levi-Civita connection; these quantities can be restricted to
slices $\RR^4_p$, $\RR^3_p$, $\MM^4_p$ to give metrics and their
Levi-Civita connections.  For any bundle $E \ra A$, we denote by
$C_A^{\infty}(E)$ the space of $C^{\infty}$ sections of $E$
defined on $A$.

\subsection{Holomorphic foliations by null planes on $\CC^4$}
\label{subsec:distnull}
By a {\em holomorphic distribution of $\alpha$-planes\/} on an
open set $A^C$ of $\CC^4$ we mean a map $\Pi$ which assigns to
each point $p$ of $A^C$ an $\alpha$-plane at $p$, $\Pi_p \subset
T_p\CC^4$\,, in a holomorphic fashion, i.e. $\Pi_p = \mbox{span}
(w_1(p), w_2(p))$ where the $w_i:A^C \ra T\CC^4$ are holomorphic.
 We can identify $\Pi$ with its image, a holomorphic subbundle of
$T\CC^4|A^C$.  

We call a holomorphic distribution $\Pi$ of $\alpha$-planes on
$A^C$ {\em autoparallel\/} if $\nabla_{w_1}w_2 \in
C^{\infty}_{A^C}(\Pi)$ for all $w_1, w_2 \in
C^{\infty}_{A^C}(\Pi)$. An autoparallel holomorphic distribution
of $\alpha$-planes has integral submanifolds which are affine
$\alpha$-planes; thus, {\em an autoparallel holomorphic
distribution of $\alpha$-planes on $A^C$ is equivalent to a
holomorphic foliation by affine $\alpha$-planes of $A^C$.} 
Similar definitions can be given replacing `$\alpha$-plane' by
`$\beta$-plane'.

\subsection{Hermitian structures on $\RR^4$} \label{subsec:distHerm}
By a {\em smooth almost Hermitian structure\/} on an open subset
$A^4$ of $\RR^4$, we mean a map $J$ which assigns to each point
$p$ of $A^4$ an almost Hermitian structure at $p$ in a smooth
fashion, i.e. $p \mapsto J_p(U_p)$ is smooth for all $U \in
C^{\infty}_{A^4}(T\RR^4)$; equivalently $J$ defines a smooth
section on $A^4$ of the bundle $E = \mbox{End}(T\RR^4) \ra
\RR^4$.  We call $J$ {\em integrable\/} if there are local
complex coordinates on $A^4$ with associated almost complex
structure $J$, this is equivalent to the vanishing of the
Nijenhuis tensor. A short calculation (see, e.g. \cite[p.
42]{Gray-Herv} or \cite[p. 169]{Salamon}) shows that this is
equivalent to
\be
\nabla^E_{JX}J =
J\nabla^E_X J \quad \forall X \in
C_{A^4}^{\infty}(T\RR^4)
\label{integ} \ee
i.e.
$$
(\nabla^E_{JX} J)(Y)
= J\big((\nabla^E_X J)(Y)\big) \quad
\forall X,Y \in C_{A^4}^{\infty}(T\RR^4)
$$
where $\nabla^E$ is the induced connection
of the bundle
$E = \mbox{End}(T\RR^4)$ given by
$(\nabla^E_X J)(Y) =
\nabla_X(JY) - J(\nabla_X Y)$. Such a $J$ is always real-analytic.

\subsection{Conformal foliations by curves on $\RR^3$}
\label{subsec:conf}
By a {\em $C^{\infty}$ (resp. $C^{\omega}$) non-zero vector
field\/} on an open subset $A^3$ of $\RR^3$ we mean a
$C^{\infty}$ (resp. $C^{\omega}$) section $U:A^3 \ra T\RR^3
\setminus \{ \mbox{zero section} \}$.  Without loss of generality
we may choose a conformal Euclidean metric $g$ on $A^3$ and scale
$U$ to be a section of the unit tangent bundle $T^1\RR^3|{A^3}$
for that metric. To such a distribution corresponds a
$C^{\infty}$ (resp. $C^{\omega}$) (oriented) foliation $\Cc$ of
$A^3$ by curves given by integrating $U$.  Note that $U$ can be
recovered from $\Cc$ as its field of (positive) unit tangents.

Let $U^{\perp}$ be the distribution of (oriented) subspaces of
$T\RR^3$ perpendicular to $U$. Then the distribution $U$ is
called {\em shear-free\/} and the corresponding foliation $\Cc$
{\em conformal\/} if Lie transport along $U$ of vectors in
$U^{\perp}$ is conformal.  In concrete terms, letting $J^{\perp}$
denote rotation through $+\pi/2$ on each oriented plane
$U_p^{\perp}, \ (p \in A^3)$, then $U$ is shear-free if and only
if ${\Ll}_U J^{\perp} = 0$ on $A^3$ where $\Ll$ denotes Lie
derivative.  Now, for any $X \in C^{\infty}(U^{\perp})$,
\begin{eqnarray*}
({\Ll}_U J^{\perp})(X)  & = & \{ ( {\Ll}_U(J^{\perp}X) \}^{\perp} -
J^{\perp} \{ {\Ll}_U(X) \}^{\perp} \\
& = & \{\nabla_U(J^{\perp}X) \}^{\perp} - \nabla_{J^{\perp}X}U -
J^{\perp} \{ \nabla_U X \}^{\perp} + J^{\perp}\nabla_X U
\end{eqnarray*}
where $\{ \ \}^{\perp}$ denotes orthogonal projection onto
$U^{\perp}$ (noting that $\nabla_X U \in C^{\infty}(U^{\perp})$
since $g(\nabla_X U, U ) = \half X g(U,U) = 0$).  Further, since
$U^{\perp}$ is a Hermitian connected bundle of rank $2$, as for
all such bundles we have
$$
\{ \nabla_U(J^{\perp}X) \}^{\perp} -J^{\perp}\{\nabla_U X
\}^{\perp} = (\nabla_U^{\mbox{\scriptsize{End}}U^{\perp}}
J^{\perp})(X) = 0 \;,
$$
hence $U$ is shear-free
if and only if
\be
\nabla_{J^{\perp}X}U = J^{\perp} \nabla_X U \;.  \label{conf}
\ee
\begin{remark} \label{re:conf-CR}
This can be interpreted as follows: A {\em CR structure\/}
\cite{Greenfield,Jacobowitz} on an odd-dimensional manifold $M =
M^{2k+1}$ is a choice of rank $k$ complex subbundle $V$ of the
complexified tangent bundle $T^cM = TM \otimes \CC$ with $V \cap
\overline{V} = \{0\}$ and
\be
[C^{\infty}(V), C^{\infty}(V)] \subset C^{\infty}(V)\;.
\label{CRint}
\ee
Given $V$ we define the Levi subbundle $H$ of $TM$ by $H \otimes
\CC = V \oplus \overline{V}$ and a real endomorphism $J:H \otimes
\CC \ra H \otimes \CC$ by multiplication by $+\ii$ (resp. $-\ii$)
on $V$ (resp. $\overline{V}$).  Conversely $(H,J)$ determines $V$
so that a CR structure can be specified by giving the pair
$(H,J)$.

Any real hypersurface $M$ of a complex manifold $\widetilde{M}$
has a canonical CR structure called the {\em hypersurface CR
structure\/} given by $V = T^{1,0}\widetilde{M} \cap T^c M$.

We give the unit tangent bundle $T^1 \RR^3 = \RR^3 \times S^2$ a
CR structure $(H,J)$ as follows: At each point $(p,U) \in \RR^3
\times S^2$ use the canonical isomorphism $\RR^3 \cong T_p\RR^3$
to regard $U^{\perp}$ as a subspace $U_p^{\perp}$ of
$T_p\RR^3$\,; then define $H_{(p,U)} = U_p^{\perp} \oplus T_U
S^2$ and define $J$ as rotation through $+\pi/2$ on $U_p^{\perp}$
together with the standard complex structure $J^{S^2}$ on $T_U
S^2$.  This is the hypersurface CR structure given by regarding
$\RR^3 \times S^2$ as a real hypersurface of the manifold ($\RR^4
\times S^2, \Jj)$ where the complex structure $\Jj$ is given by
\be
\RR^4 \times S^2 \lra
\mbox{End}\,T(\RR^4 \times S^2)\,, \qquad (p,U) \mapsto (J(U)_p,
J^{S^2})
\label{Jfund}\ee
with $J(U)_p$ the unique positive almost Hermitian structure with
$J(U)_p(\pa/\pa x_0) = U_p$ as in (\ref{JU}). A calculation using
(\ref{integ}) shows that $\Jj$ is integrable.

More generally, for any oriented Riemannian $3$-manifold $M^3$,
we can give the unit tangent bundle $T^1 M^3$ a CR structure
$(H,J)$ as follows \cite{Sato-Yam,LeBrunCR} \footnote{This is not
the CR structure on the unit tangent bundle discussed, for
example, in \cite[Chapter 7]{Blair} or \cite{Tanno}.}: At each
point $(p,U) \in T^1 M^3$ (where $p \in M^3$ and $U \in T^1_p
M^3$) the Levi-Civita connection on $M^3$ defines a splitting
$T_{(p,U)}(T^1 M^3) = T_p M^3 \oplus T_U(T^1_p M^3)$. Since
$T^1_p M^3$ is isometric to a $2$-sphere, it has a canonical
K{\"a}hler structure $J^{S^2}$.  We define $H_{(p,U)} =
U_p^{\perp} \oplus T_U (T^1_p M^3)$ and $J$ as rotation through
$+\pi/2$ on $U_p^{\perp}$ together with $J^{S^2}$ on $T_U (T^1_p
M^3)$.

For $M^3 = S^3$ this can be described more explicitly: The
differential of the canonical embedding $S^3 \hookrightarrow
\RR^4$ defines an embedding $i:T^1 S^3 \hookrightarrow T\RR^4 =
\RR^4 \times \RR^4$.  At a point $(p,U) \in T^1 S^3 \,, \ (p \in
S^3, \ U \in T^1_p S^3)$, we have
$$
di(T_{(p,U)} T^1 S^3) = \{ (X,u) \in p^{\perp} \times U^{\perp}) : \langle
X, U \rangle + \langle p, u \rangle = 0 \} \subset \RR^4 \times
\RR^4 \cong T_{(p,U)}(\RR^4 \times \RR^4) \;;
$$
then we choose $H_{(p,U)} = (U \oplus p)^{\perp} \times (U \oplus
p)^{\perp}$ and $J =$ rotation through $\pi/2$ on each plane $(U
\oplus p)^{\perp}$.  That this is integrable can be seen by
noting that any stereographic projection $S^3 \setminus
\{\mbox{point} \} \ra \RR^3$ is conformal and induces a CR
diffeomorphism between $T^1(S^3 \setminus \{\mbox{point} \}$) and
$T^1 \RR^3$.

Next give $A^3$ the CR structure defined by $(H,J) =
(U^{\perp},J^{\perp})$. Then {\em (\ref{conf}) says that $U:A^3
\ra T^1 S^3$ is CR (cf. \cite{LeBrunCR}).} \end{remark}

A concrete way of obtaining conformal foliations is given as
follows:  Let $f:A^3 \ra N^2$ be a $C^{\infty}$ (resp.
$C^{\omega}$) submersion from an open subset of $\RR^3$ to a
Riemann surface $N^2$, usually $\CC$ or $\Ci$.  Then $f$ is
called {\em horizontally conformal\/} if the restriction of its
differential $df_p$ to its {\em horizontal space\/} $(\ker
df_p)^{\perp}$ is conformal and surjective for all $p \in A^3$. 
Explicitly, $f$ is horizontally conformal if and only if, in any
local complex coordinate on $N^2$, we have on $A^3$,
\be
\left( \frac{\pa f}{\pa x_1} \right)^2 + \left(\frac{\pa f}{\pa
x_2} \right)^2 + \left( \frac{\pa f}{\pa x_3 } \right)^2
= 0\;. \label{HC}
\ee
Writing $f = f_1 + \ii f_2$ this is equivalent to
$$
g(\mbox{grad} f_1, \mbox{grad} f_2) = 0 \quad \mbox{and} \quad
g(\mbox{grad} f_1, \mbox{grad} f_1) = g(\mbox{grad} f_2,
\mbox{grad} f_2) \;.
$$ 

Then we have the simple lemma (cf. [Vai]):

\begin{lemma}
1) If $f$ is $C^{\infty}$ (resp. $C^{\omega}$) and horizontally
conformal then the foliation defined by (the fibres of) $f$ is
$C^{\infty}$ (resp. $C^{\omega}$) and conformal.

2)  All $C^{\infty}$ (resp. $C^{\omega}$) conformal foliations
are given locally in this way.
\end{lemma}

Note that the shear-free unit vector field $U$ tangent to the
conformal foliation defined by $f = f_1 + \ii f_2$ is given
by
\be
U = \mbox{grad}f_1 \times \mbox{grad} f_2 /
|\mbox{grad} f_1 \times \mbox{grad} f_2| \;.
\label{U}
\ee   

\begin{example} \label{ex:radial}
Let $f(x_1,x_2,x_3) = (x_2 \pm \ii x_3)/(|x| -x_1)$\,. It is
easily checked that $f$ is a horizontally conformal map from
$\RR^3 \setminus \{ (0,0,0) \}$ to $\Ci$.  Its level curves are
radii from the origin and give the leaves of a conformal
foliation $\Cc$ of $\RR^3 \setminus \{ (0,0,0) \}$ whose tangent
vector field is the shear-free unit vector field
\be
U(x_1, x_2, x_3) = \pm \frac{1}{\sqrt{x_1^2 + x_2^2 + x_3^2}}
(x_1, x_2, x_3) \:. \label{Urad}
\ee
Note that $U:\RR^3 \setminus \{ 0 \} \ra S^2$ is surjective.  Note
further that, when $\RR^3$ is compactified to $S^3$, this example
becomes projection from the poles: $S^3 \setminus \{(\pm
1,0,0,0)\} \ra S^2$ given by the formula $(x_0,x_1,x_2,x_3)
\mapsto \pm (x_1,x_2,x_3)/\sqrt{x_1^2 + x_2^2 + x_3^2}\;.$
\end{example}

\begin{example} \label{ex:circles}
Let $f(x_1,x_2,x_3) = -\ii x_1 \pm \sqrt{x_2^2 + x_3^2}$.
This is a horizontally conformal map from $A^3 = \RR^3 \setminus
\{ (x_1,x_2, x_3) : x_2 = x_3 = 0 \}$ to $\CC$.  Its level curves
are circles in planes parallel to the $(x_2,x_3)$-plane and
centred on points of the $x_1$-axis; these give a conformal
foliation of $A^3$ whose tangent vector field is the shear-free
unit vector field $U:A^3 \ra S^2$ given by
\be
U(x_1, x_2, x_3) = \pm \frac{1}{\sqrt{x_2^2 + x_3^2}}
(0, -x_3, x_2) \;.
\ee
Note that $U$ has $1$-dimensional image --- the equator of $S^2$.
\end{example}

\subsection{Shear-free ray congruences on $\MM^4$}
\label{subsec:SFR}
By a {\em $C^{\infty}$ (resp. $C^{\omega}$) distribution of null
directions\/} on a subset $A^M$ of Minkowski space $\MM^4$ we
mean a map $V$ which assigns to each point $p$ of $A^M$ a
$1$-dimensional null subspace of $T_p \MM^4$ in a smooth fashion,
i.e. $V_p = \mbox{span}(v(p))$ where $v:A^M \ra T\MM^4$ is a
smooth nowhere zero null vector field.  We can identify $V$ with
its image, a smooth null rank one subbundle of $T\MM^4|A^M$. 
Such a distribution integrates to a $C^{\infty}$ (resp.
$C^{\omega}$) foliation $\ell$ by null curves (i.e. curves whose
tangent vectors are all null).  We call $V$ {\em autoparallel\/}
if
\be
\nabla_v v \in C^{\infty}_{A^4}(V) \mbox{ for all } v \in
C^{\infty}_{A^4} (V) \,,
\label{aut} \ee
equivalently the integral curves of $V$ are null lines (``light
rays") and then $\ell$ is called the {\em foliation by null
lines\/} or {\em ray congruence\/}  \footnote{The term {\em ray
congruence\/} is often used to describe a three-parameter family
of rays which may not form a foliation everywhere.} on $A^M$
corresponding to $V$. Given $\ell$, the autoparallel distribution
$V$ of null directions can be recovered as its tangent field.

For a foliation $\ell$ by null lines, or the corresponding
autoparallel distribution $V$ of null directions, we define the
shear-free condition as follows:  The distribution $V^{\perp}$
orthogonal to $V$ (with respect to the conformal Minkowski metric
$g^M$) is three-dimensional and contains $V$ so that the quotient
$V^{\perp}/V$ is $2$-dimensional and $g^M$ factors to give a
positive definite inner product on it. Let $J^{\perp}$ denote
rotation through $\pi/2$ (with any orientation) on
$V^{\perp}/V\,.$  We say that $V$ is {\em shear-free} if
\be
{\cal L}_v J^{\perp} = 0
\label{SF}\ee
for any $v \in C^{\infty}(V)$\,. As in the last section this can
be interpreted as saying that Lie transport along $V$ in
$V^{\perp}/V$ is conformal.  For another interpretation, note
that (\ref{SF}) is equivalent to
\be
(\nabla_{J^{\perp}X}v)^{V^{\perp}/V} = J^{\perp}(\nabla_X
v)^{V^{\perp}/V} \mbox{  for all } X \in C^{\infty}(V^{\perp})
\label{SF2} \ee
where $( \ )^{V^{\perp}/V}$ denotes the projection $V^{\perp} \ra
V^{\perp}/V$.
More concretely, choose any complement $\Sigma$ of $V$ in
$V^{\perp}$; $\Sigma$ is called a {\em screen space}.  Then $V$
is shear-free if and only if
\be
(\nabla_{J^{\perp}X}v)^{\Sigma} = J^{\perp}(\nabla_X
v)^{\Sigma} \mbox{  for all } X \in C^{\infty}(\Sigma)
\label{SFS} \ee
where now $( \ )^{\Sigma}$ indicates projection onto $\Sigma$
along $V$.

\begin{remarks} \label{re:SFRinterp}
1. The twistorial representation of a foliation by null lines is
as follows: Recall that the affine $\alpha$-plane (\ref{I})
determined by $w = [w_0, w_1, w_2, w_3] \in \CP{3}$ intersects
the Minkowski space $\MM^4$ if and only if $w \in N^5$ (cf. Eqn.
(\ref{N5})).  Thus a foliation by null lines on an open subset
$A^M$ of $\MM^4$ is given by a map
\be
w:A^M \lra N^5 \;.  \label{w}
\ee
For a CR interpretation of condition (\ref{SFS}), see Remarks
\ref{re:SFRconf}.

2. The {\em shear tensor\/} (cf. \cite{Benn}) of a foliation
$\ell$ by null lines is defined by $S(X,Y) =$ the trace-free part
of $\half g(\nabla_X Y + \nabla_Y X, v)$ for $X,Y \in
C^{\infty}(\Sigma)$.  It is easy to see that $\ell$ is shear-free
if and only if this vanishes.  The
{\em rotation\/} or {\em twist\/} tensor is defined by $T(X,Y) =
\half g(\nabla_X Y - \nabla_Y X, v) = \half g([X,Y], v) = \half
\big(g(\nabla_Y v, X) - g(\nabla_X v, Y)\big)$.  At a point $p$ it measures
how much (infinitesimally) nearby null lines passing through the
screen space $\Sigma_p$ twist around the null line through $p$.
\end{remarks}

\subsection{Shear-free ray congruences and conformal foliations}
\label{subsec:SRFconffol}
In this section we describe a geometrical correspondence between
$C^{\infty}$ shear-free ray congruences on $\MM^4$ and
$C^{\infty}$ conformal foliations by curves on $\RR^3$.  Firstly,
recall the correspondence of Proposition \ref{prop:bijpt} at any
point $p \in \CC^4$ between unit vectors $U_p \in T_p\RR^3_p$ and
null directions $V_p$ in $T_p\MM^4_p$ given by
$$
\begin{array}{rclcl}
V_p  & = & V_p(U_p) & = & \mbox{span}(\pa/\pa t + U_p) \,, \\
U_p & = & U_p(V_p) & = & \mbox{normalized projection of } V_p
\mbox{ onto } T_p\RR^3_p\,, \\
& & & & \mbox{ i.e. the unique vector }U_p \mbox{ such that }V_p =
\mbox{span}(\pa /\pa t + U_p)\,.
\end{array}
$$
Now let $\ell$ be a $C^{\infty}$ foliation by null lines on an
open subset $A^M$ of $\MM^4$ and let $V$ be its tangent
(autoparallel) distribution of null directions as in the last
section.  Let $p \in A^M$. Then we define {\em the projection of
$\ell$ onto the slice $\RR^3_p$} as (i) the $C^{\infty}$ unit
vector field $U$ on $A^3 = A^M \cap \RR^3_p$ given by $U_q =
U_q(V_q), \ (q \in A^3)$, or, equivalently, (ii) the $C^{\infty}$
foliation $\Cc$ by curves with tangent vector field $U$.

Conversely, given a $C^{\infty}$ foliation $\Cc$ by curves or,
equivalently, a $C^{\infty}$ unit vector field $U$ on an open
subset $A^3$ of $\RR^3_p$, setting $V_q = V_q(U_q) =
\mbox{span}(\pa /\pa t + U_q)$ for $q \in A^3$ defines a
distribution of null directions on $A^3$.  We extend this to a
foliation by null lines as follows: For each $q \in A^3$, set
$\ell_q = $ the affine null line of $\MM^4$ through $q$ tangent
to $V_q\,;$ this defines a foliation $\ell$ by null lines of some
neighbourhood $A^M$ of $A^3$ in $\MM^4\,$; we define $V$ to be
its tangent distribution on $A^M$, this is an autoparallel
distribution of null directions. We call $\ell$ (or $V$) the {\em
extension} of $\Cc$ (or of $U$).  We have the following relation
between properties of $\Cc$ and $\ell$ (equivalently $U$ and
$V$):

\begin{proposition} \label{prop:SFR-conf}
Let $\Cc$ be a foliation by curves of an open subset $A^3$ of
$\RR^3_p$ and let $\ell$ be its extension to a foliation by null
lines of a neighbourhood $A^M$ of $A^3$ in $\MM^4$ as defined
above.  Then $\ell$ is shear-free on $A^M$ if and only if $\Cc$
is conformal.
\end{proposition}

\proof
We use notations as above, and write $v = \pa/\pa t + U$ so that
$V = \mbox{span}(v)$. Taking the screen space $\Sigma$ at a point
$q \in A^3$ to be $V_q^{\perp} \cap \RR^3_q = U_q^{\perp} \cap
\RR^3_q$, since $\pa/\pa t$ is parallel, the conditions
(\ref{SFS}) and (\ref{conf}) coincide so that $V$ is shear-free
at points of $A^3 = A^M \cap \RR^3_p$ if and only if $U$ is shear
free on $A^3$.  But now the Sachs' equations \cite{P-R} show that
if $V$ is shear-free on a slice cutting each null geodesic, it is
shear-free everywhere and the proposition follows.
\eproof

Recall that by a {\em germ\/} of $V$ at $A^3$ we mean an
equivalence class of distributions $V$ on open neighbourhoods of
$A^3$, deemed equivalent if they agree on a neighbourhood of
$A^3$. Then extension and projection defined above are inverse at
the level of germs, precisely:
 
\begin{theorem} \label{th:SFR-conf}
Let $p \in \MM^4$ and let $A^3$ be an open subset of $\RR^3_p$. 
Then projection onto the slice $\RR^3_p$ defines a bijective
correspondence between germs at $A^3$ of $C^{\infty}$  shear-free
ray congruences $\ell$ (equivalently shear-free autoparallel null
distributions $V$) defined on an open neighbourhood of $A^3$ in
$\MM^4$ and $C^{\infty}$ conformal foliations $\Cc$ by curves
(equivalently shear-free unit vector fields $U$) on $A^3$. The
inverse of projection is extension.
\end{theorem}

\proof
Evident from our description of the maps and the last proposition.
\eproof
     
\begin{remarks} \label{re:SFRconf}
1. Thus any $C^{\infty}$ conformal foliation $\Cc$ by curves on an
open subset of $\RR^3$ extends to a shear-free ray congruence on an open
subset $A^M$ of $\MM^4$. But then, for any point $p =
(t,x_1,x_2,x_3) \in A^M$, this projects to a $C^{\infty}$
conformal foliation by curves on an open set of the $\RR^3$-slice
$\RR^3_t \equiv \RR^3_p$ of $\MM^4$ through $p$, so that, to
$\Cc$ is associated a $1$-parameter family $\Cc_t$ of other
conformal foliations of open subsets of $\RR^3$ by curves.

Note that the transformation $\Cc \mapsto \Cc_t$
moves the vector field $U$ a distance $t$ in direction $U$ so that
the value of the unit tangent vector field $U_t$ of $\Cc_t$ at a point $q$
is equal to the value of the unit tangent vector field $U$ of $\Cc$ at
the point $q - tU_t(q)$.

2. The correspondence of the theorem can be interpreted as
follows (taking $p = 0$ for simplicity): Give $A^3$ the CR
structure $(U^{\perp}, J^{\perp})$ as
in Remark \ref{re:conf-CR}.  Then $U$ is shear-free if and only
if $U:A^3 \ra T^1 \RR^3$ is CR. Similary, by (\ref{SFS}), $\ell$ is
shear-free if and only if the map $w:A^M \ra N^5$ representing it
(Remark \ref{re:SFRinterp}) is CR when restricted to $A^3$. Now
recall that the conformal compactification $\Rt$ of $\RR^3$ may
be identified with $S^3$ and that this is contained in the
conformal compactification $\Mf$ of $\MM^4$. We have a CR
isomorphism: $k:T^1 S^3 \ra N^5$ given by sending a unit tangent
vector $U$ at a point $p$ of $S^3 = \Rt$ to the point of $\CP{3}$
representing the affine null geodesic through $p$ tangent to
$\pa/\pa t + U$.  This is given on $\RR^3$ by the map
\be
T^1\RR^3 \cong \RR^3 \times S^2 \ra N^5
\label{T1R3N5}\ee
which is the restriction of the fundamental map (\ref{fund}). 
(Since that map is holomorphic with respect to the complex
structure $\Jj$ of (\ref{Jfund}) the map (\ref{T1R3N5}), and so
also the map $k$, is CR). Then $\ell$ is represented by a map $w
= k \circ U$ so that $w$ is CR if and only if $U$ is, which
provides another proof of Proposition \ref{prop:SFR-conf}.

3. It is easily seen that $\ell$ is twist-free (Remarks
\ref{re:SFRinterp}) if and only if $\Cc$ has integrable
horizontal spaces.
\end{remarks}

\begin{example}
Consider as in Example \ref{ex:radial} the conformal foliation of
$\RR^3 \setminus \{(0,0,0)\}$ given by the shear-free vector
field
$$
U(x_1, x_2, x_3) = \pm \frac{1}{\sqrt{x_1^2 + x_2^2 +x_3^2}}(x_1,
x_2, x_3) \;.
$$
Then it is easily checked that the shear-free ray congruence
$\ell$ extending $U$ is defined on $\MM^4 \setminus
\{(t,0,0,0):t\in \RR\}$ and has tangent vector field
$v = \pa/\pa t + U$ where
\be
U(t, x_1, x_2, x_3) = \pm
\frac{1}{\sqrt{x_1^2 + x_2^2 +x_3^2}}(x_1, x_2, x_3) \;.
\label{radialSFR} \ee
Now for each $t$, \ $\ell$ projects to a conformal foliation on
$\RR^3_t$ with tangent vector field given by (\ref{radialSFR}):
note that these conformal foliations are independent of $t$.
\end{example}

\begin{example} \label{ex:circlesSFR}
Consider as in Example \ref{ex:circles} the conformal foliation
with tangent field the shear-free vector field
$$
U(x_1, x_2, x_3) = \frac{1}{\sqrt{x_2^2+x_3^2}}(0, -x_3, x_2) \;.
$$  
We compute the tangent vector field to the shear-free ray
congruence $\ell$ extending $U$. The affine null geodesic of
$\ell$ in $\MM^4$ through $(x_1, x_2, x_3)$ with direction
$\pa / \pa t + U$ is given parametrically by
\begin{eqnarray}
T  & \mapsto & \left(T, x_1, x_2 + T\left( \frac{-x_3}{\sqrt{x_2^2
+ x_3^2}} \right) , x_3 + T\left( \frac{x_2}{\sqrt{x_2^2
+ x_3^2}} \right) \right) \label{NG} \\
 & = & (T,X_1,X_2,X_3) \ , \mbox{ say.} \nonumber
\end{eqnarray}
Conversely, given $(T,X_1,X_2,X_3) \in \MM^4$ the null geodesic
of $\ell$ hits $\RR^3 \subset \MM^4$ at $(x_1, x_2,
x_3)$ where
$$
x_1 = X_1, \quad x_2 + T \left( \frac{-x_3}{\sqrt{x_2^2 + x_3^2}}
\right) =X_2, \quad x_3 + T\left( \frac{x_2}{\sqrt{x_2^2 +
x_3^2}} \right) =X_3 \;.
$$  
Solving this gives
$$
(x_1, x_2, x_3) = \frac{X_2^2+X_3^2}{R^2} \left( X_1, X_2+
\frac{T}{R}, X_3 - \frac{T}{R} \right)
$$
where $R = \sqrt{X_2^2+X_3^2-T^2}$.
Hence the tangent to the null geodesic of $\ell$ through
$(t, x_1, x_2, x_3) \in \MM^4$ is given by $v = \pa / \pa
t + U$ with
\be
U = U_t(x_1,x_2,x_3) = U(t,x_1,x_2,x_3) = \frac{r}{\sqrt{x_2^2+x_3^2}}
\left(0,-x_3 + \frac{t}{r}x_2, x_2 + \frac{t}{r}x_3 \right)
\label{E2} \ee
where $r = \sqrt{x_2^2+x_3^2-t^2}$.

Of course, again the image of $U$ is the equator of $S^2$.  For
each $t$, (\ref{E2}) gives a conformal, in fact Riemannian, foliation
$\Cc_t$ of
$\RR^3_t$.  To see what the leaf $L$ through a point $p = (x_1,x_2,x_3)$
looks like, note that, they lie in the planes $x_3 =$ constant. Fixing
$x_3$, note that, as in Remark \ref{re:SFRconf},
the value
of $U_t$ at a point $q = (x_2,x_3)$ is equal to the value of $U$ at the
point $\widetilde{q} = q-tU_t(q)$.
Since $U(\widetilde{q})$ is orthogonal to the position
vector of
$\widetilde{q}$, the points $\{0,\widetilde{q},q,q-\widetilde{q}\}$ form
a rectangle so that that the normal to $L$ at $q$ is
tangent to the circle centred at the origin with radius $|t|$. It
follows that the leaves of $\Cc_t$ are the involutes of that circle; two
such leaves are pictured in Fig. 1.
\end{example}

\subsection{Unification of the four distributions}

In \S\S \ref{subsec:distnull}, \ref{subsec:distHerm} and
\ref{subsec:SFR} we discussed the following (sets of)
distributions:

\be
\left. \begin{array}{rl}
\mbox{(D1)} & \mbox{holomorphic distributions } \Pi \mbox{ of }
\alpha \mbox{-planes on an open subset }A^C \mbox{ of }\CC^4\,,
\\

\mbox{(D2)} & \mbox{ positive almost Hermitian structures } J
 \mbox{ on an open subset } A^4 \mbox{ of } \RR^4_p\,, \\
 
\mbox{(D3)} & C^{\omega} \mbox{ distributions of null directions }
V \mbox{ on an open subset } A^M \mbox{ of } \MM^4_p\,.
\end{array} \right\}
\label{Ugerm} \ee  

(Here $p$ is a point of $\CC^4$ which we could take to be the
origin, but for later purposes we prefer not to.)

We now observe that, under complexification, these are equivalent
and that the conditions discussed in \S\S \ref{subsec:distnull},
\ref{subsec:distHerm} and \ref{subsec:SFR} then correspond.

Let $\Pi \in$ (D1).  Then, the pointwise equivalences of
Proposition \ref{prop:bijpt} give us corresponding distributions
$J = J(\Pi)$ and $V = V(\Pi)$ defined on $A^C$ by $J_q =
J_q(\Pi_q)$ etc. as in Equations (\ref{Ufor}).  Note that these
distributions are also holomorphic (with respect to the standard
holomorphic structure on $\CC^4$ given by multiplication by
$\ii \;$). They may be {\em restricted to slices\/} to give
distributions $J \in$ (D2) and $V \in$ (D3) defined on the open
sets $A^4 = A^C \cap \RR^4_p$, $A^M = A^C \cap \MM^4_p$. 
Conversely, a given $J \in$ (D2) or $V \in$ (D3) may extended to
an open subset of $\CC^4$ by analytic continuation, i.e. by
insisting that the extended quantity be holomorphic with respect
to the standard complex structure on $\CC^4$, and then $\Pi =
\Pi(J)$ or $\Pi = \Pi(V)$ is in (D1).  These operations are
clearly inverse at the level of germs.

Then we have

\begin{proposition} \label{prop:unifthree}
Let $A^C$ be a connected open subset of $\CC^4$, and let $\Pi \in$
(D1). Let $p \in A^C$ and let $A^4 = A^C \cap \RR^4_p$, $A^M =
A^C \cap \MM^4_p$. Let $J \in$ (D2), $V \in$ (D3) be defined by
restriction to slices through $p$ as above.

Then the following conditions are equivalent:

1) $\Pi$ is autoparallel on $A^C$,

2) $J$ is integrable on $A^4$,

3) $V$ is autoparallel and shear-free on $A^M$.
\end{proposition}
\proof
Firstly we show that {\em $\Pi$ is autoparallel at points of $A^4$
if and only if $J$ is integrable on $A^4$}.  To do this, let $q
\in A^4$. Then $\Pi$ is autoparallel at $q$ if and only if
$(\nabla_v w)_q \in \Pi_q$ for all $v,w \in
C^{\infty}_{A^C}(\Pi)$. But since $J$ is the $(0,1)$-subspace of
$\Pi$ this is clearly equivalent to $\nabla_v J = 0$ for all $v
\in \Pi_q$, i.e.
\be
\nabla_{X+iJX}J = 0 \quad \mbox{for all} \quad X \in T_q\RR^4_q\,.
\label{nabla} \ee
Now since $\Pi$ is holomorphic,
$$
\nabla_{iX}J = J\nabla_X J \quad \mbox{for all} \quad X \in
T_q\RR^4_q\,,
$$ 
so that (\ref{nabla}) is equivalent to
\be
\nabla_{JX}J = J \nabla_X J \quad \mbox{for all} \quad X \in
T_q\RR^4_q
\label{integ2} \ee
which is just the integrability condition (\ref{integ}) as required.

Next we show that {\em  $\Pi$ is autoparallel at points of $A^M$ if
and only if $V$ is autoparallel and shear-free on $A^M$}.  To do
this, take $q \in A^M$.  Then with $U = J(\pa/\pa x_0)$, since
$\pa/\pa x_0$ is parallel, (\ref{integ2}) is equivalent to
\be
\nabla_{JX}U = J\nabla_X U \quad \mbox{for all} \quad X \in
T_q\RR^4_q
\label{C1} \ee
and so for $v = \pa/\pa t + U$,
\be
\nabla_{X+iJX}v = 0 \quad \mbox{for all} \quad X \in
T_q\RR^4_q\;.
\label{v} \ee
Choosing $X = \pa/\pa x_0 = i\pa/\pa t$ this reads
$$
\nabla_v v = 0
$$
which is the autoparallel condition (\ref{aut}); choosing instead
$X$ in the screen space $V_q^{\perp} \cap \RR^3_q = U_q^{\perp}
\cap \RR^3_q$ gives the shear-free condition (\ref{SFS}).
Conversely, since the two choices of $X$ give vectors $X+\ii JX$
spanning $\Pi_q\,$, \big\{ (\ref{aut})  and (\ref{SFS}) \big\}
$\Rightarrow$ (\ref{v}).

To finish the proof, note that if any of the above conditions
holds at all points of $A^4$ or $A^M$, by analytic continuation
it holds throughout $A^C$.
\eproof

\begin{remark}
The proof shows that the equivalent conditions 1) to 3) of the
Proposition on quantities (D1)--(D3) are all equivalent to the
condition (\ref{C1}) on their common {\em direction vector field}
$U$ defined by (\ref{Ufor}).
\end{remark}

Consider now the following (sets of) quantities:
\be
\left. \begin{array}{rl}
\mbox{(Q1)} & \mbox{holomorphic foliations } \Ff \mbox{ by }
\alpha\mbox{-planes (equivalently, autoparallel holomorphic} \\
& \mbox{ distributions } \Pi \mbox{ of }\alpha \mbox{-planes) on an
open subset } A^C \mbox{ of } \CC^4 \,,
\\
\mbox{(Q2)} & \mbox{positive Hermitian structures } J \mbox{ on an
open subset } A^4 \mbox{ of } \RR^4_p \,,
\\
\mbox{(Q3)} & \mbox{real-analytic shear-free ray congruences } \ell
\mbox{ (equivalently, autoparallel}  \\
 & \mbox{ distributions }V \mbox{ of null directions) on an open subset }
A^M \mbox{ of }
\MM^4_p \,,
\\
\mbox{(Q4)} & \mbox{real-analytic conformal foliations } \Cc
\mbox{ by curves (equivalently, real analytic} \\
 & \mbox{ shear-free unit vector fields } U \mbox{ ) on an open subset } A^3
\mbox{ of } \RR^3_p\,.
\end{array} \right\}
\label{Q} \ee
In the last section we showed that (germs of) the quantities (Q3)
and (Q4) are equivalent in the $C^{\infty}$ category.  We now
show that all four quantities are equivalent in the $C^{\omega}$
category:

\begin{theorem} \label{th:unifgerms}
Let $\Pi \in$ (Q1).  Then restriction to slices defines
surjections (Q1) $\ra$ (Q2), (Q1) $\ra$ (Q3) and (Q1) $\ra$
(Q4).

In fact, for a fixed open set $A^3$ of $\RR^3_p$, these maps
define bijections between germs at $A^3$ of quantities (Q1),
(Q2), (Q3) and (Q4).
\end{theorem}

\proof
This follows by combining Theorems \ref{th:SFR-conf} and
\ref{prop:unifthree}. 
\eproof

Our results may be summarized by the commutative diagram where the
arrows represent restrictions to slices:

\small \be
\begin{array}{ccccc} 
 & & \vbox{ \hbox{\space \space \space \space \space \space \space
\space \space \space (Q2) $=$}
\hbox{  \{ positive Hermitian}
 \hbox{\space \space structures $J$ on $A^4$ \} } } & & \\
& \nearrow & & \searrow & \\
\vbox{ \hbox{\space \space \space \space \space \space \space \space
\space \space \space \space (Q1) $=$}
\hbox{ \{ holomorphic foliations $\Ff$}
\hbox{\space \space \space \space \space of $A^C$ by $\alpha$-planes \}
 } } & & & &
\vbox{ \hbox{\space \space \space \space \space \space \space \space \space
(Q4) $=$}
\hbox{ \{ $C^{\omega}$ foliations $\Cc$}
\hbox{\space \space of $A^3$ by curves \} } }  \\
 & \searrow & & \nearrow & \\
& &
\vbox{ \hbox{\space \space \space \space \space \space \space \space
\space \space (Q3) $=$}
\hbox{ \{ $C^{\omega}$ shear-free ray}
\hbox{congruences $\ell$ on $A^M$ \} } }     & &
\end{array}
\label{Unif} \ee \normalsize

The inverses to the maps (Q1) $\ra$ (Q2) and (Q1) $\ra$ (Q3) are
given by analytic continuation, the inverse to (Q3) $\ra$ (Q4) is
described in Theorem \ref{th:SFR-conf}, for the inverse to (Q2)
$\ra$ (Q4), see below.

This last map expresses a particular consequence of Theorem
\ref{th:unifgerms} which will be central to our work, namely that
any $C^{\omega}$ conformal foliation $U$ of an open subset of
$\RR^3$ comes from a positive Hermitian structure $J$ by {\em
projection\/}:.  Precisely:

\begin{corollary} \label{cor:JtoU}
Let $p \in \RR^4$ and let $A^3$ be an open subset of $\RR^3_p$. 
Then projection  $J \mapsto U = J(\pa/ \pa x_0)$ onto the slice
$\RR^3_p$ defines a bijective
correspondence between germs at $A^3$ of Hermitian structures $J$
defined on an open neighbourhood of $A^3$ in $\RR^4$ and
shear-free unit vector fields $U$) on $A^3$ (equivalently
$C^{\omega}$ conformal foliations $\Cc$ by curves).
\end{corollary}

\begin{remarks} \label{rem:projn}
1.  As a consequence of the theorem, given any one of the
quantities (Q1)--(Q4) we get all four distributions $\Pi \in$ (Q1),
$J \in$ (Q2), $V \in$ (Q3), $U \in$ (Q4) defined on an open subset
of $\CC^4$ and related at each point by (\ref{Ufor}) --- we call
these (or the foliations $\Ff$, $\ell$ or $\Cc$ to which they are
tangent) {\em extended quantites}.

2. The correspondence (Q2) $\ra$ (Q4) in the Corollary is the
``real analogue" of the correspondence (Q3) $\ra$ (Q4) of Theorem
\ref{th:SFR-conf}. However note that the latter applies in the
$C^{\infty}$ case too and has a inverse which is geometrically
described. To describe the inverse of  (Q2) $\ra$ (Q4) we must
follow the route (Q4) $\ra$ (Q3) $\ra$ (Q1) $\ra$ (Q2):
explicitly, given $U \in (Q4)$, extend $J = J(U)$ to an open
subset of $\MM^4$ by insisting that it be constant along the null
lines of the shear-free null geodesic congruence $\ell$ which
extends $U$ (\S \ref{subsec:SRFconffol}), then extend to an open
subset of $\CC^4$ by analytic continuation (with respect to the
standard complex structure on $\CC^4$ given by $\ii$), then
restrict to $\RR^4_p$.

3. The restriction $\Ff \in $ (Q1) $\mapsto \ell \in $ (Q3) can
be described geometrically: the null lines of $\ell$ are the
intersections of the $\alpha$-planes of $\Ff$ with $\MM^4_p$.  Note
however, that only those $\alpha$-planes represented by points of
$N^5_p$ (see \S \ref{subsec:twistors}) have non-empty intersection
with $\MM^4_p$ so
that the {\em inverse\/} (Q3) $\ra$ (Q1) is not completely
described in this way.

4. The inverse of the map (Q1) $\ra$ (Q2), $\Ff \mapsto J$ can be
described in a purely geometrical way similar to the inverse of
(Q3) $\ra$ (Q4) (see \S \ref{subsec:SRFconffol}): given $J \in$
(Q2), define $\Pi = \Pi(J)$ on $A^4$ by (\ref{Ufor}).  Then the
null planes of $\Ff$ are those tangent to the distribution
$\Pi$.

5. Given a {\em real-analytic} conformal foliation $\Ff$ by curves
of an open subset $A^R$ of $\RR^3$ by curves, there corresponds a
$5$-parameter family of real-analytic conformal foliations by
curves of open subsets of $\RR^3$ got by extending $\Ff$ to a
holomorphic foliation by $\alpha$-planes and then projecting this
to $\RR^3$-slices.
\end{remarks}

\subsection{Other formulations of the equations and a Kerr theorem}
\label{subsec:Kerr}
Given any quantity (Q1)--(Q4) and so all four extended quantities $\Pi$,
$J$, $V$, $U$ related by (\ref{Ufor}) writing $U = \sigma^{-1}(\ii \mu)$
where $\mu
= w_1/w_0 \in \Ci$, we see that (\ref{C1}) is equivalent to
the equation
\be
\nabla_Z \mu = 0 \quad \mbox{for all} \quad Z \in \Pi_p \;.
\label{nablamu} \ee
Recall that, in the coordinates $(z_1, \widetilde{z}_1, z_2,
\widetilde{z}_2)$ (see (\ref{coords})),  a basis for $\Pi$ =
is given by (\ref{NT0}) or
(\ref{NT}).  Since $\Pi$ is the $(0,1)$-tangent space of $J = J(\Pi)$
it follows that $(z_1 - \mu \widetilde{z}_2, z_2 +
\mu \widetilde{z}_1)$ are holomorphic with respect to $J$ and, so,
provide local complex coordinates for $J$. Furthermore
(\ref{nablamu}) reads
\be
\left( \frac{\pa}{\pa \widetilde{z}_1} - \mu \frac{\pa}{\pa z_2}
\right) \mu = 0 \ ,
\quad
\left( \frac{\pa}{\pa \widetilde{z}_2} + \mu \frac{\pa}{\pa z_1}
\right) \mu = 0 \;.
\label{EC} \ee
These equations restrict on real slices to the equations
\cite[Eqns. (6.1)]{Wood-4d} or \cite[Eqns. (2.1)]{Baird-JMP}:
\be
\left( \frac{\pa}{\pa \overline{z}_1} - \mu \frac{\pa}{\pa z_2}
\right)\mu = 0 \;,
\quad
\left( \frac{\pa}{\pa \overline{z}_2} + \mu  \frac{\pa}{\pa z_1}
\right)\mu = 0 \;,
\label{ER} \ee
and on Minkowski slices, writing $v' = x_1 + t$, $v = x_1 - t$, to
\be
\left( \frac{\pa}{\pa v} - \mu \frac{\pa}{\pa z_2} \right)\mu = 0
\ , \quad
\left( \frac{\pa}{\pa \overline{z}_2} + \mu \frac{\pa}{\pa v'}
\right) \mu = 0 \;,
\label{EM} \ee
which are the equations for a shear-free ray congruence as given
by \cite{Penrose}, \cite[II (7.4.6)]{P-R} or \cite[p. 50]{H-T}
modulo conventions.

Now recall (\S \ref{subsec:twistors}) that the set of affine
$\alpha$-planes in $\Cf$ can be identified with $\CP{3}$ so
that a distribution $\Pi$ of $\alpha$-planes may be represented
by a map $w:A^C \ra \CP{3}$ with $w(p) \in \widehat{p}$ for all
$p \in A^C$. Then it is not hard to see that (\ref{C1}) is
equivalent to
\be
\mbox{holomorphic part of } dw(Z) = 0 \mbox{ for all } p \in
A^C \mbox{ and } Z \in \Pi_p \;.
\label{dw} \ee
This formulation needs no coordinates and makes sense on the
compactified space $\Cf$ hence {\em Proposition
\ref{prop:unifthree}, Theorem \ref{th:unifgerms} and Corollary
\ref{cor:JtoU} extend to compactified spaces}. We thus have a
unified {\em Kerr Theorem\/} for our quantities:

\begin{proposition}
Given any of the quantities (Q1)--(Q4), representing its extension by
a holomorphic map $w:A^C \ra \CP{3}$ on an open subset of $\CC^4$
with $w(p) \in \widehat{p}$ for all $p \in A$, there is a unique
complex hypersurface $S$ of $\CP{3}$ which is the image of $w$. 

Conversely, given a complex hypersurface $S$ of $\CP{3}$, any
holomorphic map $w:A^C \ra \CP{3}$ with $w(p) \in \widehat{p}$
and image in $S$ defines related (extended) quantities
(Q1)--(Q4).

Further, if $S$ is given by
\be
\psi(w_0, w_1, w_2, w_3) = 0
\label{S} \ee
where $\psi$ is homogeneneous and holomorphic, then the direction
vector field of the quantity away from points at infinity is
given by $U = \sigma^{-1} (\ii w_1/w_0)$ where $[w_0, w_1] = \mu
(z_1, \widetilde{z}_1, z_2, \widetilde{z}_2)$ is a solution to
\be
\psi (w_0, w_1, w_0 z_1 - w_1 \widetilde{z}_2, w_0 z_2 + w_1
\widetilde{z}_1) = 0 \;.
\label{K} \ee   
\end{proposition}
\proof
Given $w$, set $S =$ the image of $w$. That $S$ is a complex
hypersurface follows from (\ref{dw}). Indeed, in local
coordinates, if $[w_0,w_1] = \mu(z_1, \widetilde{z}_1, z_2,
\widetilde{z}_2)$ represents the direction field of the quantity,
then by (\ref{EC}), $[w_0,w_1]$ is a holomorphic function of
$(z_1 - \mu \widetilde{z}_2, z_2 + \mu \widetilde{z}_1)$. The
latter pair equals $(w_2,w_3)$ by the incidence relations
(\ref{I}). Hence $[w_0,w_1]$ are holomorphic functions of
$[w_2,w_3]$ and an equation of the form (\ref{S}) is satisfied
showing explicitly that $S$ is a complex hypersurface. 
Conversely, if $S$ is a complex hypersurface, then it is locally
of the form (\ref{S}) and reversing the above argument shows that
$\mu = w_1/w_0$ satisfies (\ref{EC}) and so $U = \sigma^{-1}(\ii
\mu) \in$ (Q4) giving related quantities $\Pi \in$ (Q1), $J \in$ (Q2)
and $V \in$ (Q3).
\eproof

The surface $S$ is called the {\em twistor surface\/} of the
quantity.

\begin{remark} \label{re:bij-twistor}
The bijections (Q3) or (Q4) to (Q1) or (Q2) can be described
twistorially as follows (taking $p = 0$ for simplicity):
By  \ref{re:SFRconf}, \ $\ell \in$
(Q3) or $\Cc \in$ (Q4) defines a CR map with image a
$3$-dimensional CR submanifold $N^3$ of $N^5$; this is
\cite{PenCR} the intersection of a complex hypersurface $S$ with
$N^5$; $S$ then determines $\Ff \in (Q1)$ and $J \in (Q2)$ as
above.

Note that if $\ell$ or $\Cc$ is only $C^{\infty}$ but has
non-integrable horizontal spaces everywhere, then it can be
extended to a Hermitian structure $J$ on one side of $\RR^3$;
precisely, let $\RR^4_+$ (resp. $\RR^4_-)$ be the set
$\{ (x_0,x_1,x_2,x_3)
\in \RR^4 : x_0 > \mbox{( resp.} < 0 ) \}$, then there exists a
$C^{\infty}$ Hermitian structure on an open subset $A^4$ of
$\RR^4_+ \cup \RR^3$ or $\RR^4_- \cup \RR^3$ with $A^4 \cap \RR^3
= A^3$ which is $C^{\omega}$ on $A^4 \setminus A^3$.  To see this
note that the CR manifold $N^3 \subset \CP{3}$ has non-zero Levi
form \cite[II, p.220 ff.]{P-R} and so, by a theorem of Harvey and
Lawson \cite{Har-Law}, is the boundary of a complex hypersurface
$S$.  This means that $J$ can be defined on one side of $\RR^3$.
\end{remark}


\section{Harmonic morphisms, null planes and SFR congruences}
In this section we define Minkowski and complex- harmonic
morphisms and show that, except in trivial cases, their equations
can be reduced to the equations (\ref{EM}), (\ref{EC}) for
shear-free ray (SFR) congruences or holomorphic foliations by
$\alpha$-planes
\subsection{Harmonic morphisms on Euclidean and Minkowski space}
For any (semi-)Riemannian manifolds $(M^m,g)$, $(N^n,h)$, a {\em
harmonic morphism} $\phi:M^m \ra N^n$ is a map which pulls back
germs of harmonic functions on $N^n$ to germs of harmonic
functions on $M^m$.  By \cite{Fu1,Ish} for the Riemannian case
and \cite{Fu2} for the semi-Riemannian case, these can be
characterized as harmonic maps $\phi$ which are {\em horizontally
weakly conformal}, i.e. at each point $p \in M$ either the
differential $d\phi_p = 0$ or the restriction of $d\phi_p$ to the
horizontal space $(\ker d\phi_p)^{\perp}$ is conformal and
surjective. In particular, a smooth map $\phi:A^4 \ra \CC$ from
an open subset of Euclidean space $\RR^4$ is a harmonic morphism
if it satisfies Laplace's equation
$$
\Delta \phi \equiv \frac{\pa^2\phi}{\pa x_0^2} + \frac{\pa^2\phi}{\pa
x_1^2} + \frac{\pa^2\phi}{\pa x_2^2} + \frac{\pa^2\phi}{\pa x_3^2} = 0
$$
and the horizontal weak conformality condition:
\be
\left( \frac{\pa\phi}{\pa x_0}\right)^2 + \left(
\frac{\pa\phi}{\pa x_1} \right)^2 + \left( \frac{\pa\phi}{\pa
x_2}\right)^2 + \left( \frac{\pa\phi}{\pa x_3}\right)^2 = 0
\label{HWCEucl} \ee
whereas a smooth map $\phi:A^M \ra \CC$ from an open subset of
Minkowski space $\MM^4$ is a {\em Minkowski harmonic morphism\/},
i.e. a harmonic morphism with respect to the Minkowski metric, if
and only if it satisfies the {\em wave equation}:
\be
\Box\phi \equiv - \frac{\pa^2\phi}{\pa t^2} + \frac{\pa^2\phi}{\pa
x_1^2} + \frac{\pa^2\phi}{\pa x_2^2} + \frac{\pa^2\phi}{\pa
x_3^2} = 0
\label{wave} \ee
and the horizontal weak conformality condition:
\be
- \left( \frac{\pa\phi}{\pa t}\right)^2 + \left(
\frac{\pa\phi}{\pa x_1} \right)^2 + \left( \frac{\pa\phi}{\pa
x_2}\right)^2 + \left( \frac{\pa\phi}{\pa x_3}\right)^2 = 0
\label{HWCMink} \ee
at all points $(t,x_1,x_2,x_3) \in \MM^4$, i.e. {\em a Minkowski
harmonic morphism is a ``null" solution to the wave equation}
(cf. \cite{Bawospin}).  Note that both (\ref{HWCEucl}) and
(\ref{HWCMink}) express horizontal weak conformality with respect
to any metric conformally equivalent to the standard metrics.

For a harmonic morphism between (semi-)Riemannian manifolds there
are three types of point $p$:

a) $d\phi_p = 0$; we call such points {\em critical points};

b) $d\phi_p \neq 0$ and $\ker d\phi_p$ is non-degenerate;

c) $d\phi_p \neq 0$ and $\ker d\phi_p$ is degenerate.  Then \cite{Fu2}
$(\ker d\phi_p)^{\perp} \subset \ker d\phi_p\,$.

If the last case occurs at some point $p$ we call $\phi$ {\em
degenerate (at $p$)}; a simple example of a harmonic morphism
$\phi:\MM^4 \ra \CC$ which is degenerate everywhere is given by
$\phi(t,x_1,x_2,x_3) = t-x_1$ or, more generally, $f(t-x_1)$ for
any smooth function $f:\RR \ra \CC$.  Note that this has
$1$-dimensional image.  For more theory and examples, see, for
example  \cite{Baird-Luminy,Wood-Sendai} for the Riemannian case,
and \cite{Parmar,Fu2} in the semi-Riemannian case. Note, in
particular, that, in the Riemannian case, a smooth harmonic
morphism is always real-analytic \cite{E-S}, but the above
example shows that this is not so in the semi-Riemannian case.

The equations for a harmonic morphism to a $2$-dimensional
codomain are conformally invariant in the codomain, i.e. if
$\phi:M^m \ra N^2$ is a harmonic morphism to a $2$-dimensional
Riemannian manifold and $\rho:N^2 \ra N'^2$ is a weakly conformal
map to another $2$-dimensional Riemannian manifold, the
composition $\rho \circ \phi$ is a harmonic morphism. Thus the
pair of equations (\ref{wave},\ref{HWCMink}) makes sense for a
map $\phi:\MM^4 \ra N^2$ to a Riemann surface.

\subsection{Complex-harmonic morphisms}
By a {\em complex-harmonic function\/} $\phi:A^C \ra \CC$ on an
open subset $A^C$ of $\CC^4$ we mean a {\em holomorphic map\/}
satisfying the complex Laplace's equation:
\be
\sum_{i=0}^3 \frac{\pa^2\phi}{\pa x_i^2}  =  0\,, \quad 
 (x_0,x_1,x_2,x_3) \in A^C.   \label{cxha}
\ee 

By a {\em complex-harmonic morphism} $\phi:A^C \ra \CC$ we mean a
holomorphic map satisfying (\ref{cxha}) and a complex version of
the horizontal weak conformality condition:
\be
 \sum_{i=0}^3 \left( \frac{\pa \phi}{\pa x_i} \right)^2  =  0\,, \quad
 (x_0,x_1,x_2,x_3) \in A^C.
\label{HWCcx} \ee

We may adapt the arguments of \cite{Fu1,Ish} to characterize
complex-harmonic morphisms $A^c \ra \CC$ as those holomorphic
maps which pull back holomorphic functions to complex-harmonic
functions. Note also that, since these
equations are conformally invariant in $\phi$, we can replace
$\CC$ by any Riemann surface.

\begin{proposition} \label{prop:cxha}
a) Let $\phi:A^C \ra N^2$ be a complex-harmonic morphism from an
open subset $A^C \subset \CC^4$ to a Riemann surface.  Then, for
any $p \in A^C$,

1) $\phi|{A^C\cap \RR^4_p}$ is a harmonic morphism (w.r.t. the
Euclidean metric);

2) $\phi|{A^C\cap \MM^4_p}$ is a harmonic morphism (w.r.t the
Minkowski metric),

b) All harmonic morphisms from open subsets of $\RR^4$ to Riemann
surfaces and real-analytic harmonic morphisms from open subsets
of $\MM^4$ to Riemann surfaces arise in this way.
\end{proposition} \proof

a) Immediate from the equations.

b) As noted above, any harmonic morphism from an open subset of
$\RR^4$ to a Riemann surface is $C^{\omega}$.  By analytic
continuation, this is the restriction of a holomorphic map on an
open subset of $\CC^4$ and this holomorphic map is
complex-harmonic.  The $\MM^4$ case is similar.
\eproof

An important example of a harmonic morphism is the following:

\begin{proposition} \label{prop:E-hamorph} Let $\mu:A \ra \Ci$
where $A$ is an open subset of $\CC^4$  (resp. $\RR^4$, $\MM^4$)
satisfy equations (\ref{EC}) (resp. (\ref{ER}), (\ref{EM})). 
Then $\mu$ is a complex-harmonic morphism (resp. harmonic
morphism, Minkowski harmonic morphism).
\end{proposition}
\proof
We give the proof for the complex case only. To establish
horizontal weak conformality (\ref{HWCcx}) we have, from
(\ref{EC}),
$$
\frac{\pa\mu}{\pa\widetilde{z}_1}\frac{\pa\mu}{\pa z_1}
+ \frac{\pa\mu}{\pa\widetilde{z}_2}\frac{\pa\mu}{\pa z_2} =
\mu\frac{\pa\mu}{\pa z_2}\frac{\pa\mu}{\pa z_1}
-\mu\frac{\pa\mu}{\pa z_1}\frac{\pa\mu}{\pa z_2} = 0 \ ;
$$
complex-harmonicity (\ref{cxha}) follows from
$$
\frac{\pa^2\mu}{\pa z_1 \pa\widetilde{z}_1}
+ \frac{\pa^2\mu}{\pa z_1 \pa\widetilde{z}_1}
= \frac{\pa}{\pa z_1} \left( \mu \frac{\pa\mu}{\pa z_2} \right)
- \frac{\pa}{\pa z_2} \left( \mu \frac{\pa\mu}{\pa z_1} \right)
= 0 \;.
$$

Thus the direction field $U$ (or its representative $\mu = -\ii
\sigma(U)$) of any of our quantities (Q1) -- (Q4) in (\ref{Q}) gives
a harmonic morphism; further examples are provided by postcomposition of
one of these with a holomorphic or antiholomorphic map
to a Riemann surface.  We shall see that, away from
critical points, this gives all non-trivial harmonic morphisms
locally.  Firstly, recall

\begin{theorem} \label{th:Wood-4d} Let $\phi:A^R
\ra N^2$ be a harmonic morphism without critical points
from an open subset $A^R$ of
$\RR^4$ to a Riemann surface.  Then there
exists a Hermitian structure on $A^R$ such that $J$ is parallel
along each connected component of a fibre of $\phi$. Further, for
any $p \in A^R\,$, there is a neighbourhood $A_1^R$ of $p$ in $A^R$,
and a holomorphic map $\rho:V \ra \Ci$ from an open subset $V$ of
$N^2$ such that $\mu = \rho \circ \phi$ represents (\S
\ref{subsec:repcoord}) $J$.
\end{theorem}

We have the following analogue for complex-harmonic morphisms:

\begin{theorem} \label{th:cxhamorphSFR} Let $\phi:A^C \ra N^2$ be
a complex-harmonic morphism without critical points
from an open subset $A^C$ of $\CC^4$
to a Riemann surface.  Then there exists
a holomorphic foliation $\cal{F}$ of $A^C$ by $\alpha$-planes or
by $\beta$-planes such that each connected component of the
fibres of $\phi$ is the union of parallel null planes of $\cal
F$. Further, for any $p \in A^C$ there is a neighbourhood $A_1^C$
of $p$ in $A^C$, and a holomorphic map $\rho:V \ra \Ci$ from an
open subset $V$ of $N^2$ such that $\mu = \rho \circ \phi$
represents (\S \ref{subsec:repcoord}) the direction field of
$\Ff$.
\end{theorem}

\proof
Let $p \in A^C$.  By Proposition \ref{prop:cxha},  \ $\phi$
restricts to a harmonic morphism on $A^4 = A^C \cap \RR^4_p$
which is easily seen to be submersive.  By \cite{Wood-4d},
$\phi|{A^4}$ is holomorphic with respect to some Hermitian
structure $J$ which is constant along each connected component of
a fibre of $\phi$.  Replacing $z_2$ by $\widetilde{z}_2$ if
necessary, we can assume that $J$ is positively oriented.
Representing $J$ by $\mu:A^4 \ra \Ci$, we have
\be
\frac{\pa\phi}{\pa \widetilde{z}_1} - \mu \frac{\pa\phi}{\pa z_2} = 0 \,,
\quad
\frac{\pa\phi}{\pa \widetilde{z}_2} + \mu \frac{\pa\phi}{\pa z_1} = 0 
\,,\label{muphi}\ee
\be
\frac{\pa\mu}{\pa \widetilde{z}_1} - \mu \frac{\pa\mu}{\pa z_2} = 0 \,,
\quad
\frac{\pa\mu}{\pa \widetilde{z}_2} + \mu \frac{\pa\mu}{\pa z_1} = 0  
\,, \label{ER2} \ee
at points of $A^4$, the first equation expressing holomorphicity
of $\phi$ with respect to $J$ (cf. \cite{Wood-4d}) and the
second, integrability of $J$ (cf. (\ref{ER})). Extend $\mu$ to
$A^C$ by (\ref{muphi}), noting that this well-defines $\mu$ since
not all the partial derivatives of $\phi$ can vanish
simultaneously: $\mu$ is then a holomorphic function.  By
analytic continuation, (\ref{ER2}) holds at all points of $A^C$
so that $\mu$ defines a holomorphic foliation $\Ff$ by
$\alpha$-planes.  By (\ref{muphi}) $\phi$ is constant along any
$\alpha$-plane of $\Ff$ so that each fibre of $\phi$ is the union
of $\alpha$-planes of $\Ff$.  Further, $J$ and so $\mu|{A^4}$ is
constant along the connected components of fibres of
$\phi|{A^4}$; by analytic continuation, $\mu:A^C \ra \Ci$ is
constant along the connected components of the fibres of
$\phi:A^C \ra \Ci$ so that the $\alpha$-planes of $\Ff$ making up
a connected component of a fibre of $\phi$ are all parallel.

Lastly, since $\mu$ is constant on the leaves of the foliation
given by the fibres of $\phi$, it factors through local leaf
spaces as $\mu = \rho \circ \phi$.  Since $\phi$ and $\mu$ are
both holomorphic (with respect to $\ii$ and $J$), $\rho$ must be
holomorphic. \eproof

We have an analogue for Minkowski harmonic morphisms:

\begin{theorem} \label{th:SFR}
Let $\phi:A^M \ra N^2$ be a real-analytic Minkowski harmonic
morphism without critical points
from an open subset $A^M$ of $\MM^4$ to a Riemann
surface.  Then there is a shear-free ray
congruence $\ell$ on $A^M$ such that each connected component of
a fibre of $\phi$ is the union of parallel null lines of $\ell$.
Further, for any $p \in A^c$ there is a neighbourhood $A_1^M$ of
$p$ in $A^M$, and a holomorphic map $\rho:V \ra \Ci$ from an open
subset $V$ of $N^2$ such that $\mu = \rho \circ \phi$ represents
the direction field of $\ell$.
\end{theorem}

\proof
Extend $\phi$ to an open subset $A^C$ of $\CC^4$ with $A^C \cap
\MM^4 = A^M$.  Let $\Ff$ be the holomorphic foliation by null
planes of Theorem \ref{th:cxhamorphSFR}.  For each $q \in A^M$,
set $\ell_q = $ the intersection of the leaf of $\Ff$ through $q$
with $\MM^4$.  By Theorem \ref{prop:SFR-conf}, $\ell$ is
shear-free and is the desired ray congruence.  Lastly, $\rho$ is
given as before.
\eproof

Note that in the real and complex cases the condition ``without
critical points" is equivalent to ``$\phi$ is submersive". This is not so
in the Minkowski case where $\phi$ may be degenerate. We can actually
be more precise in this case:

\begin{corollary} \label{cor:deg}
Let $\phi:A^M \ra N$ be a real-analytic harmonic morphism from an
open subset $A^M$ of $\MM^4$.  Suppose that $\phi$ is degenerate
at $p$ with $d\phi_p \neq 0$.  Then there exists a unique null
direction $V_p \in T_p\MM^4$ such that $V_p \subset \ker
d\phi_p$. Furthermore, $\ker d\phi_p = V_p^{\perp}$.  If,
further, at each point $q$ in the connected component of the
fibre through $p$, \ $\phi$ is degenerate with $d\phi_q \neq 0$,
then that connected component is the affine null 3-space tangent
to $V_p^{\perp}$.
\end{corollary}

\proof
By \cite{Fu2}, $(\ker d\phi_p)^{\perp} \subset \ker d\phi_p$. This
means that $\ker d\phi_p$ must be three-dimensional.  But then
$(\ker d\phi_p)^{\perp}$ is $1$-dimensional and null, say $(\ker
d\phi_p)^{\perp} = V_p$.  So $\ker d\phi_p = V^{\perp}_p$.

To prove uniqueness of $V_p$, suppose that $V'_p \subset \ker d\phi_p$
is another null direction. Then $V'_p \subset V_p^{\perp}$ which
is easily seen to imply that $V'_p =V_p$.

This means that the distribution $p \mapsto V_p$ must be tangent
to the SFR congruence of the theorem, and so each $V_p$ is
parallel for all $p$ in a connected component of a fibre. The
last assertion follows from the fact that the connected
component of the fibre is $3$-dimensional and has every tangent
space parallel to $V_p$.
\eproof

\begin{remark} \label{re:totgeod}
If $\phi$ is non-degenerate at $p$, the kernel of $d\phi_p$ has
Minkowski signature and so contains two null lines.  If the fibre
at $p$ is totally geodesic, then it is the affine Minkowski plane
tangent to that kernel. If all the fibres are totally geodesic,
then the two null lines both give SFR congruences $\ell$ as in
the theorem; if, however, at least one fibre is not totally
geodesic, the SFR congruence $\ell$ of the theorem is unique. 
Similar remarks hold in the $\RR^4$ and $\MM^4$ cases.
\end{remark}

Call a smooth map from an open subset $A^4$ of $\RR^4$ to a
Riemann surface {\em K\"ahler\/} if it is holomorphic with
respect to a K\"ahler structure on $A^4$, such maps are always
harmonic morphisms \cite{Fu1}.  Analogously, call a smooth map $\phi$ from
an open subset $A^C$ of $\CC^4$ (resp. $\MM^4$) to a Riemann
surface {\em K\"ahler\/} if the restriction of $\phi$ (resp. of the
analytic continuation of $\phi$) to real slices is K\"ahler.
Again such maps are always complex- (resp. Minkowki) harmonic
morphisms. Clearly, all such maps are explicitly known.  If the
$\mu$ in any of the last three theorems is constant, then $\phi$
is K\"ahler (the converse is true unless $\phi$ has totally
geodesic fibres in which case it can be holomorphic w.r.t. a
K\"ahler structure of one orientation and also w.r.t. a non-K\"ahler
Hermitian
structure of the other orientation, Example \ref{ex:Minhamorphs}
(4) is of this type).  Then we have a converse to Proposition
\ref{prop:E-hamorph}:

\begin{corollary} \label{cor:compo}
Let $\phi:A \ra N^2$ be a non-K\"ahler harmonic (resp.
complex-harmonic, real-analytic Minkowski harmonic) morphism from
an open subset of $\RR^4$ (resp. $\CC^4$, $\MM^4$) to a Riemann
surface.  Then, for any $p \in A$ with $d\phi_p \neq 0$,
there is a neighbourhood $A_1$ of $p$ and a
non-constant holomorphic map $\rho:V \ra \Ci$ from an open subset
of $N^2$ such that $\mu = \rho \circ \phi$ satisfies Equations
(\ref{ER}) (resp. (\ref{EC}), (\ref{EM})).
\end{corollary}
Thus, except in the trivial K\"ahler case, the harmonic morphism equations can be reduced to the first order equations (\ref{ER}), (\ref{EC}) or (\ref{EM}).  
\begin{examples} \label{ex:Minhamorphs}

1) The simplest examples of Minkowski harmonic morphisms from
$\MM^4$ to $\CC$ are given by a) $(t,x_1,x_2,x_3) \mapsto x_2
+\ii x_3$ which is non-degenerate everywhere and surjective and
b) $(t,x_1,x_2,x_3) \mapsto x_1 - t$ which is degenerate
everywhere and has $1$-dimensional image $\RR$. Note that in both
cases, the fibres are totally geodesic; in a) the SFR congruences
of Theorem \ref{th:SFR} (two by Remark \ref{re:totgeod}) have
leaves with (null) directions $(1, \pm 1,0,0)$, in
b) the SFR congruence (just one by Corollary \ref{cor:deg}) has
leaves with (null) direction $(1,1,0,0)$.

2) Consider the harmonic morphism $\phi:\MM^4 \ra \RR^2$ given by
$(t,x_1,x_2,x_3) \mapsto (x_1-t+x_2^2-x_3^2,2x_2 x_3)$.  For any
real $c$, the fibre $\phi= c$ is the union of two surfaces
intersecting on the line $x_1 = t, \ x_2 = x_3 = 0$; that line
consists of degenerate points; other points on the fibre are
non-degenerate.

3) The direction field $U$ of the SFR extending the simple
conformal foliation of Example \ref{ex:circles} is given by (\ref{E2}). It
defines a Minkowski harmonic morphism from the cone $A^M = \{
(t,x_1,x_2,x_3) \in \MM^4: x_2^2+x_3^2 > t^2 \}$ to $S^2$ which
is degenerate everywhere and has image the equator of $S^2$. 
Note that $U = \sigma (\ii\mu)$ with
$$
\mu(t,x_1,x_2,x_3) =
\frac{r}{x_2^2+x_3^2} \left\{ \left( x_2 +
\frac{t}{r}x_3 \right) + \ii \left(x_3 - \frac{t}{r}x_2 \right)
\right\} $$ where $r = \sqrt{x_2^2+x_3^2-t^2}$ so that $\mu$ is a
harmonic morphism $A^M \ra \CC$, also degenerate, with image the
unit circle. The fibre of $U$ or $\mu$ through any point $p$ is
the affine plane perpendicular to $U_p$ tangent to $U_p$,
$\pa/\pa x_1$ and the vector in the $(x_2,x_3)$-plane
perpendicular to $U_p$.

4) Define $\mu:\MM^4 \setminus \{(t,x_1,x_2,x_3): x_2 = x_3 = x_1
+ t = 0 \} \ra \Ci$ by $$ \mu = \frac{x_2 + \ii x_3}{i(x_1 + t)}
\;. $$ This defines an SFR and so a Minkowski harmonic morphism.
It is non-degenerate except at points on the 3-plane
$\{(t,x_1,x_2,x_3): x_1 + t = 0 \} .$

5) Take the twistor surface $S = \{[w_0,w_1,w_2,w_3] \in \CP{3} :
w_0 w_1 + w_2 w_3 = 0 \}$.  Then $\mu$ satisfies $$ \mu + (z_1 -
\mu \widetilde{z}_2)(z_2 + \mu \widetilde{z}_1) = 0 $$ and any
local smooth solution gives the direction field of an SFR and
thus a Minkowski harmonic morphism on an open subset of $\MM^4$. 
It is clear that this harmonic morphism does not have totally
geodesic fibres.

Note that in Examples 3) to 5) the $\rho$ of Theorem \ref{th:SFR}
is the identity map whereas for Examples 1) and 2) it is the
constant map these examples being K\"ahler. \end{examples}


\section{Finding conformal foliations} \label{sec:theory}
\subsection{The boundary of a hyperbolic harmonic morphism}
According to Corollary \ref{cor:JtoU}, given a positive Hermitian
structure $J$ on an open subset $A^4$ of $\RR^4$, if $U$ is its
direction vector field, i.e. $U = J(\pa/\pa x_0)$, then, for any
$p \in A^4$, \ $U|{A^3}$ is the tangent vector field of a
$C^{\omega}$ conformal foliation on the open set $A^3 = A^4 \cap
\RR^3_p$ and all $C^{\omega}$ conformal foliations by curves of
open subsets of $\RR^3$ arise this way.  To get an explicit
description of all conformal foliations we need to find the
horizontally conformal functions $f$ on $\RR^3$ whose level sets
give the leaves of the foliation. To this end note the

\begin{proposition} \label{prop:restr-R3}
(i) Let $\phi:A^4 \ra \CC$ be a submersive map from an open subset
of $\RR^4$ which is holomorphic with respect to a positive
Hermitian structure $J$ on $A^4$. Let $p \in A^4$, \ $U =
J(\pa/\pa x_0)$ be the correponding shear-free vector field on
the open subset $A^3 = A^4 \cap \RR^3_p$ of $\RR^3_p$ and $\Cc$
the conformal foliation of $A^3$ given by its integral curves. 
If
\be
\frac{\pa \phi}{\pa x_0} = 0 \mbox{ on } A^3
\label{orthog} \ee
then $f= \phi|A^3$ is a real-analytic horizontally conformal
submersion which is constant on the leaves of $\Cc$.

(ii) All real-analytic horizontally conformal submersions (and
so, all real-analytic conformal foliations by curves) on open
subsets of $\RR^3$ are given in this way.
\end{proposition}
\proof
(i) Let $q \in A^3$.  By holomorphicity, since $U_q = J_q (\pa/\pa
x_0)$, \ the directional derivative $U_q(f) = 0$ so that $f$ is
constant on the leaves of $\Cc$.  If $\{ e_2, e_3 = Je_2 \}$ is a
basis for $U_q^{\perp} \cap \RR_p^3\,$, holomorphicity of $\phi$
implies that $e_3(f) = \ii e_2(f)$ so that $f$ is horizontally
conformal.  Submersivity of $f$ easily follows from that of
$\phi$.

(ii) Given a real-analytic horizontally conformal submersion
$f:A^3 \ra \CC$ on an open subset of $\RR^3$, let $U$ be the unit
vector field tangent to its level curves. By Corollary
\ref{cor:JtoU}, there is an open subset $A^4$ of $\RR^4$ with
$A^3 = A^4 \cap \RR^3$ and a positive Hermitian structure $J$ on
$A^4$ such that $U = J(\pa / \pa x_0)$ on $A^3$. Then $f$ is CR
with respect to the hypersurface structure on $A^3$ induced by
$J$ and so may be extended to a function $\phi$ holomorphic with
respect to $J$ on a subset $A_1^4$ of $A^4$ with $A^3 = A_1^4
\cap \RR^3$. Then, since at points of $A^3$, \ $U = J(\pa / \pa
x_0)$ and $U(\phi) = 0$, we have $\pa\phi / \pa x_0 = 0$ on $A^3$
establishing (ii).
\eproof

\begin{remark} \label{rem:extn}
The extension of $f$ to $\phi$ can be described more geometrically
in the same way as the extension of $U$ to $J$ (see Remarks
\ref{rem:projn}), namely: Extend $f$ to an open subset of $\MM^4$
by insisting that it be constant along the null lines of the
shear-free null geodesic congruence $\ell$ which extends $U$ (see
\S \ref{subsec:SRFconffol}), then extend to an open subset of
$\CC^4$ by analytic continuation (with respect to the standard
complex structure on $\CC^4$ given by multiplication by $\ii$)
and restrict to $\RR^4$. The extension of $f$ can also be
described in a twistorial way as follows: As in Remark
\ref{re:bij-twistor}, $U$ defines a $3$-dimensional CR
submanifold $N^3$ of $N^5$ which is the intersection of $N^5$ and
a complex hypersurface $S$ of $\CP{3}$ --- the twistor surface of
$J$.  Then $f$ defines a CR function $F = f \circ \pi$ on $N^3$
which we may extend to a holomorphic function $\Phi$ on a
neighbourhood of $N^3$ in $S$; on a possibly smaller
neighbourhood, $\Phi$ is of the form $\phi \circ \pi$ for a
unique function $\phi$ on an open subset of $\RR^4$.  This is the
desired extension.
\end{remark}       

We now show how to find such functions $\phi$ using harmonic
morphisms. Equip $\breve{\RR}^4 \equiv \RR^4 \setminus \RR^3$
with the {\em hyperbolic metric} $g^H = \left( \sum_{i=0}^{3}
dx_i^2 \right) / x_0^2$ so that each component $\RR_+^4$,
$\RR_-^4$ is isometric to hyperbolic $4$-space $\HH^4$. Let
$\breve{A}^4$ be an open subset of $\breve{\RR}^4$, then we call
a smooth map $\phi:\breve{A}^4 \ra \CC$  a {\em hyperbolic
harmonic map} if it is a harmonic map with respect to the
hyperbolic metric $g^H$.  This holds if and only if
\be
x_0 \sum_{i=0}^4 \frac{\pa^2\phi}{\pa x_i^{{}2}} -
2 \frac{\pa\phi}{\pa x_0} = 0
\label{hyp-ha} \ee
at all points of $\breve{A}^4$.

Similarly, $\pi:\breve{A}^4 \ra \CC$ will be called a {\em
hyperbolic harmonic morphism\/} if it is a harmonic morphism with
respect to the metric $g^H$, such maps are characterized as
satisfying (\ref{hyp-ha}) and (\ref{HWCEucl}).

Recall (\cite{Wood-4d}, see also \cite{Baird-JMP})
\begin{theorem} \label{th:Wd-4d2} (i) Any submersive hyperbolic
harmonic morphism $\phi:\breve{A^4} \ra \CC$ is holomorphic with
respect to some Hermitian structure $J$ on $\breve{A^4}$ and has
superminimal fibres with respect to $J$, i.e. $\ker d\phi \subset
\ker \nabla^H J$ on $\breve{A}^4$  where $\nabla^H$ is the
Levi-Civita connection of the hyperbolic metric on
$\breve{\RR}^4$.

(ii) Conversely, let $J$ be a Hermitian structure on an open
subset $\breve{A}^4$ of $\breve{\RR^4}$, and $\phi:\breve{A^4}
\ra \CC$ a non-constant map which is holomorphic with respect to
$J$, then $\phi$ is hyperbolic harmonic if and only if, at points
where $d\phi \neq 0$, its fibres are superminimal with respect to
$J$.
\end{theorem}

To formulate this analytically, let $\Theta$ be the homogeneous
holomorphic contact form
\be
\Theta = w_1 dw_2 - w_2 dw_1 - w_0 dw_3 + w_3 dw_0
\label{HCF} \ee
on $\CP{3}$. Then $\ker\Theta$ gives the horizontal spaces of the
restriction of the twistor map (\ref{twistormap}):  $\pi:\CP{3}
\setminus N^5 \ra (\breve{\RR}^4,g^H)$. Set $\Phi = \phi \circ
\pi$. Then $\phi$ is a hyperbolic harmonic morphism if and only
if
\be
\ker d\Phi \subset \ker \Theta\,.
\label{superminimal} \ee
Equivalently, let $w:\breve{A}^4 \ra \CP{3}$ be the section of
$\pi$ corresponding to $J$, i.e. with image the twistor surface
$S$ of $J$ (see \S \ref{subsec:Kerr});  then we may pull back
$\Theta$ to a $1$-form $\theta = w^*(\Theta)$ on $A^4$. 
Condition (\ref{superminimal}) now reads $\ker d\phi \subset \ker
\theta$,  this condition ensuring that the fibres of $\phi$ are
superminimal.

Note that if $A^4$ is an open subset of $\RR^4$ or of $\RR_+^4
\cup \RR^3$ (or $\RR_-^4 \cup \RR^3$) with $A^4 \cup \RR^3$
non-empty, any ($C^2$, say) map which is a hyperbolic harmonic
map on $\breve{A}^4 = A^4 \setminus \RR^3$ satisfies
(\ref{hyp-ha}) and so (\ref{orthog}) at the boundary $A^3 = A^4
\cup \RR^3$. A key property for us is the following converse:

\begin{proposition} \label{prop:orthog-hypha}
Let $A^4$ be a connected open subset of $\RR^4$, $\RR_+^4 \cup
\RR^3$ or $\RR_-^4 \cup \RR^3$ with $A^3 = A^4 \cap \RR^3$
non-empty, and let $\phi:A^4 \ra \CC$ be a non-constant $C^1$ map
which is holomorphic with respect to a Hermitian structure $J$ on
$\breve{A}^4 = A^4 \setminus \RR^3$ and submersive at almost all
points of $A^3$. Then $\phi$ satisfies (\ref{orthog}) on $A^3$ if
and only if $\phi|\breve{A}^4$ is a hyperbolic harmonic
morphism.
\end{proposition}

\proof
It suffices to work at points where $\phi$ is submersive. At such
points note that (\ref{orthog}) holds if and only if $\ker \phi =
\mbox{span} (\pa/\pa x_0 , J \pa/\pa x_0)$. Now let $S$ be the
twistor surface of $J$ and let $\Phi:S \ra \CC$ be defined by
$\Phi = \phi \circ \pi$.

We have the
\begin{lemma} 
The pull-back $\theta = w^* \Theta$ to $A^4$ satisfies $\mbox{\rm
span} (\pa/\pa x_0 , J \pa/\pa x_0) \subset \ker \theta$ at all
points of $A^3$.
\end{lemma}
\proof
Since the (complexified) normal to $\RR^3$ is given by the
annihilator of $\mbox{span} \{dz_1 - d\widetilde{z}_1,
dz_2,d\widetilde{z}_2 \}$, it suffices to show that on $\RR^3 =
\{ z_1 + \widetilde{z}_1 = 0 \}$,\ $\theta$ is a linear
combination of those three forms.  To do this, taking
differentials in (\ref{I}) gives \begin{eqnarray*} dw_2 & = & z_1
dw_0 + w_0 dz_1 - \widetilde{z}_2 dw_1 - w_1 d\widetilde{z}_2
\;,\\ dw_3 & = & z_1 dw_0 + w_0 dz_1 - \widetilde{z}_2 dw_1 - w_1
d\widetilde{z}_2 \;.
\end{eqnarray*}
Substituting these into (\ref{HCF}) and rearranging gives
$$
\theta = (w_3 - w_0 z_2 + w_1 z_1)dw_0 - (w_2 + w_0 \widetilde{z}_1
+w_1 \widetilde{z}_2)dw_1 + w_0 w_1 (dz_1 - d\widetilde{z}_1) -
w_0^{{}2} d\widetilde{z}_1 -
w_1^{{}2} d\widetilde{z}_2 \;.
$$
But, by (\ref{I}), the coefficients of $dw_0$ and $dw_1$ vanish
when $z_1 + \widetilde{z}_1 = 0$ and the lemma follows. \eproof

By the lemma the condition (\ref{orthog}) is equivalent to the
superminimality of the fibres of $\phi$ at points of $A^3$, viz.:
$\ker d\phi \subset \ker \theta$, or, equivalently $\ker d\Phi
\subset \ker \Theta$ on the real hypersurface $N^3 = w(A^3)$ of
$S$.  But this is a holomorphic condition, so by analytic
continuation, if $\phi$  has superminimal fibres at points of
$A^3$ then it has superminimal fibres on the whole of $A^4$ and
we are done.
\eproof

Combining Propositions \ref{prop:restr-R3} and
\ref{prop:orthog-hypha} we obtain the basis of our method:

\begin{theorem} \label{th:HWC-hyp}
Let $f:A^3 \ra \CC$ be a real-analytic horizontally conformal
submersion on an open subset of $\RR^3$. Then there is an open
subset $A^4$ of $\RR^4$ with $A^4 \cap \RR^3 = A^3$ and a
real-analytic submersion $\phi:A^4 \ra \CC$  with $\phi|A^3 = f$
such that $\phi|A^4 \setminus \RR^3$ is a hyperbolic harmonic
morphism. In fact $\phi \mapsto f = \phi|A^3$ defines a bijective
correspondence between germs at $A^3$ of real-analytic
submersions $\phi:A^4 \ra \CC$ on open neighbourhoods of $A^3$ in
$\RR^4$ which are hyperbolic harmonic on $A^4\setminus \RR^3$ and
real-analytic horizontally conformal submersions $f:A^3 \ra
\CC$.
\end{theorem}

\begin{remarks}

1. The hyperbolic harmonic morphism $\phi$ has totally geodesic
fibres if and only if the level sets of $f$ are circles, see
\cite{Baird-circles}.

2. Since the map $\phi$ is found geometrically (Remark
\ref{rem:extn}), there is no problem with K\"ahler points as
there was in \cite{Wood-4d,Baird-JMP}.

3. If $f$ is only $C^{\infty}$, then, as in Remark
\ref{rem:projn}, if the distribution $\mbox{span} \{ \mbox{grad}
f_1, \mbox{grad} f_2 \}$ is nowhere integrable, we can extend $f$
to one side of $\RR^3$; precisely, there is an open subset $A^4$
of $\RR^4_+ \cup \RR^3$ or $\RR^4_- \cup \RR^3$ with $A^4 \cap
\RR^3 = A^3$ and a $C^{\infty}$ map $\phi:A^4 \ra \CC$  with
$\phi|A^3 = f$ such that $\phi|A^4 \setminus \RR^3$ is a
hyperbolic harmonic morphism.
\end{remarks}  

\begin{corollary}
Let $c$ be a embedded real-analytic curve in $\RR^3$.  Then there
is an embedded real-analytic surface $s$ in an open subset $A^4$
of $\RR^4$ which is minimal on $A^4 \setminus \RR^3$ with respect
to the hyperbolic metric, hits $\RR^3$ orthogonally and has $s
\cap \RR^3 = c$.  In fact there is a conformal real-analytic
foliation of an open subset $A^4$ of $\RR^4$ by such surfaces
with $s$ a leaf.
\end{corollary}

\proof
Firstly, $c$ can be embedded in a conformal foliation by curves of
an open subset of $\RR^3$ as follows:  construct the normal
planes to $c$ and integrate the vector field given by the normals
to these.  This gives a foliation on an open neighbourhood of $c$
in $\RR^3$ which has totally geodesic integrable horizontal
spaces and so (see, for example, \cite{Wood-contemp}) is
Riemannian.  (To get a conformal foliation which is not
Riemannian, replace the planes by spheres, possibly of varying
radii.) Representing this foliation as the level curves of a
real-analytic horizontally conformal submersion $f:A^3 \ra \CC$
on an open subset of $\RR^3$, construct a hyperbolic harmonic
morphism $\phi$ as in the theorem: its fibres give the desired
real-analytic foliation.
\eproof

\begin{remark}
It is easily seen from the equations that any $C^1$ surface in an
open subset $A^4$ of $\RR^4$, $\RR_+^4 \cup \RR^3$ or $\RR_-^4
\cup \RR^3$ which is minimal with respect to the hyperbolic
metric on $\breve{A^4} = A^4 \setminus \RR^3$ hits $\RR^3$
orthogonally.
\end{remark}

\subsection{Finding horizontally conformal functions from twistor
surfaces} \label{subsec:findingHWC}

We now show how to use Theorem \ref{th:HWC-hyp} to find explicitly
the horizontally conformal submersion and conformal foliation by
curves on an open subset of $\RR^3$ corresponding to a given
complex hypersurface $S$ of $\CP{3}$.  In fact, by introducing a
parameter $a \in \CC^4$, we can obtain a $5$-parameter family of
horizontally conformal submersions whose level sets give the
$5$-parameter family of conformal foliations of open subsets of
$\RR^3$ of in Remarks \ref{rem:projn}.  For this we need to
translate the hyperbolic metric to different slices:  First
recall (\S \ref{subsec:twistors}) the map $\pi_a:\CP{3} \ra
\Rf_a$ defined by $w \mapsto$ the intersection of the
$\alpha$-plane (\ref{I}) determined by $w$ with $\Rf_a$.  Then,
denoting the first component of $a$ by $a_0$, we have

\begin{lemma}
Let $a \in \CC^4$. Equip $\breve{\RR}^4_a = \RR^4_a \setminus
\RR^3_a$ with the hyperbolic metric $g_a^H = \left( \sum_{i=0}^3
dx_i^{{}2} \right) /(x_0- \Re a_0)^2$ and let $N^5_a =
\pi_a^{-1}(\RR^3_a)$.  Then the kernel of the holomorphic contact
form \be \Theta_a = \Theta_{a_0} = -2a_0(w_1 dw_0 - w_0 dw_1) +
w_1 dw_2 - w_2 dw_1 - w_0 dw_3 + w_3 dw_0 \label{Theta_a} \ee
restricted to $\CP{3} \setminus N^5_a$ gives the horizontal
distribution of $\pi_a:\CP{3} \setminus N^5_a \ra
(\breve{\RR}_a^4, g^H_a)$.
\end{lemma}
\proof
Writing $a$ in coordinates (\ref{coords}) as
$(a_1,\widetilde{a_1}, a_2, \widetilde{a_2})$, it can easily be
checked that translation $T_a:x \mapsto x+a$ in $\CC^4$
corresponds to the map $\widetilde{T}_a: \CP{3} \ra \CP{3}$ given
by $[w_0,w_1,w_2,w_3] \mapsto [w_0,w_1,w_2 + a_1 z_0 -
\widetilde{a}_2 z_1, z_3 + a_2 z_0 + \widetilde{a}_1 z_1]\,$,
that is, $\pi_a \circ \widetilde{T}_a = T_a \circ \pi$. Then
$\Theta_a = (\widetilde{T}_a^{-1})^*\Theta = \widetilde{T}_{-a}^*
\Theta$; on calculating this, (\ref{Theta_a}) follows.
\eproof

Now let $S \subset \CP{3}$ be a given complex hypersurface and $a
\in \CC^4$. Let $U$ be an open set in $S$ such that $\pi_a$ maps
$U$ diffeomorphically onto an open set $A^4$ of $\RR_a^4$. Then
$U$ defines a Hermitian structure $J$ on $A^4$ represented by the
section $w:A^4 \ra U$ of $\pi_a$. If $S$ is given by
\be
\psi(w_0,w_1,w_2,w_3) = 0
\label{psi} \ee
(where $\psi$ is homogeneous holomorphic) $w$ is given by solving
(\ref{K}).  Given a holomorphic map $\zeta:U \ra \Ci$, set
$\phi_a = \zeta \circ w:A^4 \ra \Ci$. Then $\phi_a$ is
holomorphic with respect to $J$ and so, by Theorem
\ref{th:Wd-4d2}, is a harmonic morphism with respect to the
hyperbolic metric on $\breve{\RR}_a^4$ if and only if the level
surfaces of $\zeta$ are horizontal, i.e. tangent to
$\ker\Theta_a$.

Set $\Sigma^S_a = \{w \in S : \ker\Theta_a = TS \}$.  Then on $S
\setminus \Sigma^S_a$, \ $\ker\Theta_a \cap TS$ is a
one-dimensional holomorphic distribution so that its integral
(complex) curves foliate $S \setminus \Sigma^S_a$. Note that
$\pi_a(\Sigma^S_a) \cap A^4$ is the set of K\"ahler points of $J$
\ $= \{ p \in A^4 : \nabla^H_v J = 0 \ \forall v \in T_p A^4 \}$
where $\nabla^H$ is the Levi-Civita connection of the hyperbolic
metric on $\breve{\RR}^4_a$.

To find $\phi_a$ we proceed as follows: Firstly let $c:V \ra
\CC^2$, \ $[w_0,w_1,w_2,w_3] \mapsto (\zeta,\eta)$ be complex
coordinates for $S$ on an open subset $V$ of $U$. Then we can
solve (\ref{I}) locally to find the composition $c \circ w$, \ $p
\mapsto (\zeta(p),\eta(p))$.

Next let $\widetilde{\zeta} = \widetilde{\zeta}(\zeta,\eta)$ be a
holomorphic function with $\pa\widetilde{\zeta}/ \pa\zeta \neq 0$
and set $\widetilde{\eta} = \eta$.  Then $(\zeta, \eta) \mapsto
(\widetilde{\zeta}, \widetilde{\eta})$ is locally a complex
analytic diffeomorphism.  Then, $\phi_a(p) =
\widetilde{\zeta}(\zeta(p),\eta(p))$ restricts to a hyperbolic
harmonic morphism on $\breve{\RR}^4_a$ if and only if the level
sets of $\widetilde{\zeta}$ are superminimal, the condition for
this is
\be
\Theta_a \left( \frac{\pa}{\pa \widetilde{\eta}} \right) = 0 \,.
\label{h_a} \ee
By the chain rule,
$$
\frac{\pa}{\pa \widetilde{\eta}} = \frac{\pa}{\pa \eta} +
\frac{\pa \zeta}{\pa\widetilde{\eta}}
\frac{\pa}{\pa \zeta} \quad \mbox{and} \quad \frac{\pa \zeta}{\pa
\widetilde{\eta}} = -
\frac{\pa\widetilde{\zeta}}{\pa \eta}
\Big/ \frac{\pa\widetilde{\zeta}}{\pa \zeta}
$$
so that (\ref{h_a}) reads
\be
\Theta_a \left( \frac{\pa}{\pa \zeta} \right)
\frac{\pa\widetilde{\zeta}}{\pa \eta}
- \Theta_a \left( \frac{\pa}{\pa \eta} \right)
\frac{\pa\widetilde{\zeta}}{\pa \zeta} =
0 \;. \label{E} \ee
This equation can be solved to get a holomorphic function
$\widetilde{\zeta} = \widetilde{\zeta}_a (\zeta,\eta)$ such that
$\phi_a:p \mapsto \widetilde{\zeta}_a (\zeta(p),\eta(p))$
restricts to a hyperbolic harmonic morphism on $(\breve{\RR}^4_a,
g^H_a)$.  Note that the solution is not unique, but, as the level
sets of $\widetilde{\zeta}_a$ are uniquely determined on $V$, \
$\phi_a$ is unique up to postcomposition with a holomorphic
function.  Restriction of $\phi_a$ to $\RR^3_a$ then gives a
horizontally conformal submersion $f$ on an open subset of
$\RR^3_a$ and hence the conformal foliation by curves
corresponding to $S$ and $a$. Note that the method actually finds
a holomorphic function $\phi_a$ on an open subset of $\CC^4$,

In summary, given a complex hypersurface $S$ of $\CP{3}$ there is
defined locally a real Hermitian structure $J$ and all the
related quantities (\ref{Q}), in particular, $U = J(\pa /\pa
x_0)$. Then given any $a \in \CC^4$, the integral curves of $U$
defined a $C^{\omega}$ conformal foliation on an open subset of
the slice $\RR^3_a$ and all such foliations are given this way.
Then we have shown above how to find a holomorphic function
$\phi_a:A^C \ra \CC$ on an open subset of $\CC^4$ such that the
level curves of $\phi_a|\RR^3_a$ give the leaves of the conformal
foliation.


\section{Examples} \label{sec:examples}
\subsection{Linear Examples}

Let $S$ be the $\CP{2}$ given by the homogeneous linear function
$$
\psi(w_0,w_1,w_2,w_3) = b_0 w_0 + b_1 w_1 + b_2 w_2 +b_3 w_3
$$
where $(b_0,b_1,b_2,b_3) \neq (0,0,0,0)$.  
We consider a representative example where $\psi$ is given by
$$
[b_0,b_1,b_2,b_3] = [0,s,0,1]
$$
with $s$ a real parameter.  Note that $s = 0$ if and only if
$[b_0,b_1,b_2,b_3] \in N^5$, this will be a special case.
Equation (\ref{K}) reads
$$
s \mu + (z_2 + \mu \widetilde{z}_1) = 0
$$
so that the direction field $U$ of the quantities
corresponding to $S$ is given by $U = \sigma^{-1}(\ii\mu)$ where
$$
\mu = -z_2/(\widetilde{z}_1 + s)\,.
$$
Parametrize $S\setminus \{z_0 = 0 \}$ by
$$
(\zeta, \eta) \mapsto [1,w_1,w_2,w_3] =[1,\zeta,\eta, -s\zeta] \,.
$$
Equation (\ref{HCF}) reads
\begin{eqnarray*}
\theta_a  & = & 2a_0 dw_1 + w_1 dw_2 - w_2 dw_1 - dw_3 \\
   & = & 2a_0 d\zeta + \zeta d\eta - \eta d\zeta + sd\zeta
\end{eqnarray*}
so that Equation (\ref{E}) reads
$$
 (-2a_0 + \eta - s)\frac{\pa\widetilde{\zeta}}{\pa\eta} +
\zeta \frac{\pa \widetilde{\zeta}}{\pa\zeta} = 0
$$
which has a solution
$$
\widetilde{\zeta} = - (-2a_0 + \eta - s)/\zeta \;.
$$
The incidence relations (\ref{I}) read
\be \left. \begin{array}{rcl}
z_1 - \zeta \widetilde{z}_2 & = & \eta \\
z_2 + \zeta \widetilde{z}_1 & = & -s\zeta \end{array} \right\}
\ee
with solution
$$
\zeta   =  -\frac{z_2}{\widetilde{z}_1 + s} \,, \qquad
\eta  =  \frac{z_1 \widetilde{z}_1 + z_2 \widetilde{z}_2 +
z_1 s}{\widetilde{z}_1 + s}
$$
so that
$$
\phi_a = - \frac{ z_1 \widetilde{z}_1 + z_2 \widetilde{z}_2
-2a_0(\widetilde{z}_1 + s ) + (z_1 -\widetilde{z}_1)s -s^2}{z_2}
\;.$$
Putting $a_0 = -\ii c$ gives the horizontally conformal function
on the slice $\RR^3_c = \RR^3_{a_0}$:
$$
\phi_c = - \frac{(x_1+c)^2 + x_2^2 + x_3^2 - s^2
     + 2\ii(x_1 + c)s}{x_2+\ii x_3} \;.
$$
Writing $\rho = |s|$, if $s \geq 0$, this is the composition of
$$
\RR^3  \stackrel{T^1_c}{\lra}  \RR^3 
\stackrel{\sigma_3^{-1}}{\lra} S^3(\rho)
\stackrel{\overline{H}}{\lra}  S^2  \stackrel{\sigma}{\lra}  \Ci
$$
where $T^1_c$ is the translation $(x_1,x_2,x_3) \mapsto
(x_1+c,x_2,x_3)$, \ $\sigma_3^{-1}$ is the inverse of
stereographic projection from $(-\rho,0,0,0)$ given by
$(X_1,X_2,X_3) \mapsto \big((\rho^2 - |X|^2, 2\rho X_1, 2\rho
X_2, 2\rho X_3\big)/(\rho^2 + |X|^2) $ and $\overline{H}$ is the
(``conjugate") Hopf map given by $\CC^2 \supset S^3(\rho) \ni
(z_1,z_2) \mapsto -\overline{z}_1/z_2$.  If $s < 0$ we replace
$H$ by the Hopf map $(z_1,z_2) \mapsto -z_1/z_2$.  In either case
the fibres of $\phi$ give the conformal foliation of $\RR^3$ by
the circles of Villarceau, see \cite{Wilker} for a description,
\cite[II p. 62]{P-R} or the cover of {\em Twistor Newsletter\/}
for a picture, and compare with a different treatment in
\cite{Baird-circles}. For $s = 0$ the foliation degenerates to a
bunch of circles tangent to the $x_1$ axis at $(-c,0,0)$.

We remark that an arbitrary linear function $\psi$ gives a
foliation by circles but not all foliations by circles are given
this way.

\subsection{Quadratic examples}
We now give three examples where the twistor surface is quadratic
so that we can solve all equations explicitly.
\begin{example}
Let the twistor surface be the quadratic surface $S =
\{ [w_0,w_1,w_2,w_3] \in \CP{3} : w_0 w_3 - w_1 w_2 = 0 \}$. Then
the direction field $U$ of the corresponding quantities is given by
$U = \sigma^{-1} (\ii \mu)$ where
$$
z_2 + \mu \widetilde{z}_1 - \mu (z_1 - \mu \widetilde{z}_2) = 0
$$
which has solutions
\be
\mu  =  \frac{z_1 - \widetilde{z}_1 \mp
\sqrt{(z_1-\widetilde{z}_1)^2
- 4z_2 \widetilde{z}_2}} {2\widetilde{z}_2}
\label{mu1} \ee
so that
$$
\ii \mu   = 
\frac{- x_1 \pm |x|}{x_2 - \ii x_3}
 =  \frac{x_2 + \ii x_3}{x_1 \pm |x|}
$$
where $|x| = \sqrt{x_1^2+x_2^2+x_3^2} = \sigma^{-1}(\ii\mu)\,.$
Thus on any slice  $\RR^3_t$, \ $U$ is given by
$$
U(t,x_1,x_2,x_3) = \pm (x_1,x_2,x_3)
 / |x|
$$
which gives the tangent vector field of foliation by radial lines of
Example \ref{ex:radial}.

Carrying out the calculations of \S \ref{subsec:findingHWC} gives
$\phi_a = \mu$ for all $a$ --- we omit the details which are
similar to the next example. Then, for each $a$, \ $\sigma^{-1}
\circ \phi_a$ restricted to $\breve{\RR}^4_a$ gives the harmonic
morphism $(\breve{\RR}^4_a  \setminus \{x_0 \mbox{-axis} \},
g_a^H ) \ra S^2$ given by orthogonal projection to $\RR^3_a$
followed by radial projection.
\end{example}

\begin{example}
This time let the twistor surface be the quadratic surface $S = \{
[w_0,w_1,w_2,w_3] \in \CP{3} : w_0 w_3 + w_1 w_2 = 0 \}$. Then
the direction field $U$ of the corresponding quantities is given
by $U = \sigma^{-1}(\ii \mu)$ where
$$
z_2 + \mu \widetilde{z}_1 + \mu (z_1 - \mu \widetilde{z}_2) = 0
$$
which has solutions
$$
\mu  =  \frac{z_1 + \widetilde{z}_1 \pm
\sqrt{(z_1+\widetilde{z}_1)^2 + 4z_2 \widetilde{z}_2}}
{2\widetilde{z}_2} =  \frac{x_0 \pm s}{x_2 - \ii x_3}
$$
where $s = \sqrt{x_0^2 + x_2^2 + x_3^2}\,.$ \
Note that
$$
\mu|\RR^3_0 = \pm \frac{\sqrt{x_2^{{}2} +x_3^{{}2}}}{x_2-\ii x_3}
= \pm \frac{x_2+\ii x_3}{\sqrt{x_2^{{}2} +x_3^{{}2}}}
$$
so that on $\RR^3_0$,
$$
U(x_1,x_2,x_3) = \pm \frac{1}{\sqrt{x_2^{{}2} +x_3^{{}2}}}
(0,-x_3,x_2)
$$
which gives the tangent vector field of the conformal foliation
$\Cc$  by circles around the $x_1$-axis of Example
\ref{ex:circles}. Further note that $\mu|\MM^4$ is given by
\begin{eqnarray*}
\mu(t,x_1,x_2,x_3) & = & \frac{-\ii t \pm r}{x_2-\ii x_3} \\
& = & \frac{r}{x_2^{{}2} + x_3^{{}2}}
\left\{ \left( \pm x_2 + \frac{t}{r}x_3 \right) + \ii
\left(\mp x_3 -\frac{t}{r}x_2 \right) \right\}
\end{eqnarray*}
where $r = \sqrt{x_2^2 +x_3^2 - t^2}$
and so, on the open set $x_2^2 +x_3^2 > t^2$, we have
$$
U(t,x_1,x_2,x_3) = \frac{r}{\sqrt{x_2^{{}2} + x_3^{{}2}}} \left(
0, \pm x_3 + \frac{t}{r}x_2, \pm x_2 + \frac{t}{r} x_3 \right)
\,;
$$
then $v = U + \pa/\pa t$ gives the tangent field of the shear-free
congruence $\ell$ extending $\Cc$ discussed in Example
\ref{ex:circlesSFR}.

Parametrizing $S$ away from $w_0 = 0$ by $(\zeta,\eta) \mapsto
[1,w_1,w_2,w_3] = [1,\eta,-\zeta, \zeta \eta]$, the incidence
relations (\ref{I}) read
$$
\left. \begin{array}{ccr}
z_1 - \eta \widetilde{z}_2 & = & - \zeta \\
z_2 + \eta \widetilde{z}_1 & = & \zeta \eta
\end{array} \right\}
$$
which have solutions
$$
\zeta = -\ii x_1 \pm s \quad \mbox{and} \quad
\eta = \frac{x_0 \pm s}{x_2 - \ii x_3} \,. 
$$  
On $\MM^4$,
\be
\zeta = -\ii x_1 \pm r \quad \mbox{and} \quad \eta =
\frac{-\ii t \pm r}{x_2-\ii x_3}
\label{zeta-eta-1} \;.\ee
These solutions are real analytic except on the cone $x_2^2 +
x_3^2 = t^2$.  We obtain smooth solutions if we avoid this
cone. Note that $S \cap N^5$ has equation:
$(\zeta+\overline{\zeta})(|\eta|^2-1) = 0$ and so consists of the
union of two submanifolds:  $S_1:\Re \zeta = 0$ and $S_2: |\eta|
= 1$.  Smooth branches of $(\zeta,\eta)$ will correspond to CR
maps $A \ra S \cap N^5$ with image lying in $S_1 \setminus S_2$
or $S_2 \setminus S_1$; \ $S_1 \cap S_2$ corresponding to the
branching set of $(\zeta,\eta)$ which is the subset of $A$ on
which the square root $r$ vanishes.  Specifically, note from
(\ref{zeta-eta-1}) that if $t > |z_2|$, \ $\Re \zeta = 0$
corresponding to points of $S_1$ and if $t < |z_2|$, \ $|\eta| =
1$ corresponding to points of $S_2$.

We have
\begin{eqnarray*}
\theta_a & = & 2a_0 dw_1 + w_1 dw_2 - w_2 dw_1 - dw_3 \\
& & 2a_0 d\eta - \eta d\zeta + \zeta d\eta - \zeta d\eta+\eta d\zeta
\end{eqnarray*}
so that (\ref{E}) reads
$$
\eta \frac{\pa \widetilde{\zeta}}{\pa \eta} + a_0
\frac{\pa\widetilde{\zeta}}{\pa \zeta} = 0 $$
which has a solution
\be
\widetilde{\zeta}_a = \zeta - a_0 \ln \eta \;.
\label{sol1} \ee
Setting $\phi_a = \widetilde{\zeta}_a \big(\zeta(p),\eta(p) \big)$
gives the complex-valued map on a dense subset of $\CC^4$:
\be 
\phi_a(x_0,x_1,x_2,x_3) = -\ii x_1 \pm s
- a_0 \ln \frac{x_0 \pm s}{x_2-\ii x_3} \,.
\label{phia2} \ee
For any $a \in \CC^4$ this restricts to a complex-valued harmonic
morphism on a dense subset of $(\breve{R}_a^4, g_a^H)$. In
particular, when $a = 0$, this simplifies to the hyperbolic
harmonic morphism on $\breve{R}^4\setminus \{x_1 \mbox{-axis}\}$
given by
\be
\phi(x_0,x_1,x_2,x_3) = -\ii x_1 \pm \sqrt{x_0^2+x_2^2+x_3^2}
\label{circles-hamorph} \ee
which further restricts on $\RR^3$ to
\be
\phi(x_1,x_2,x_3) = -\ii x_1 \pm \sqrt{x_2^{2} + x_3^{2}} \;,
\label{L} \ee
the level surfaces of this giving the conformal foliation by
circles round the $x_1$-axis of  Example \ref{ex:circles}. The
harmonic morphism (\ref{circles-hamorph}) has fibres given by the
Euclidean spheres having these circles as diameters, these
spheres are totally geodesic in $(\breve{R}^4, g^H)$.

For definiteness, now take plus signs in the above.
Then, putting $a_0 = -\ii t$ in (\ref{phia2}) and restricting to the
open set $\{(x_1,x_2,x_3) : x_1^2 + x_2^2 > t^2 \} \subset \RR^3_t
= \RR^3_a$  gives the horizontally conformal map
\begin{eqnarray*}
\phi_t  =  \phi_a & = & -\ii x_1 + r + \ii t \ln
\frac{r-\ii t}{x_2- \ii x_3} \\
       & = & -\ii x_1 + r - t \arg \frac{r- \ii t}{x_2- \ii x_3} \;,
\end{eqnarray*}
the level curves of this lie in the planes $x_1 =$ constant \ and are
the involutes of the unit circle pictured in Fig. 1.  
\end{example}

\begin{example}
This time, as in Example \ref{ex:Minhamorphs} let the twistor
surface be the quadratic surface $S = \{ [w_0,w_1,w_3,w_4] \in
\CP{3} : w_0 w_1 + w_2 w_3 = 0 \}$.  Then $\mu$ satisfies
$$
\mu + (z_1 - \mu \widetilde{z}_2)(z_2 + \mu \widetilde{z}_1) = 0
$$
so that
$$
\mu = \frac{(1 + z_1 \widetilde{z}_1 - z_2 \widetilde{z}_2) \pm s}
{2\widetilde{z}_1\widetilde{z}_2}
$$
where $s = \sqrt{(1 - z_1 \widetilde{z}_1 - z_2 \widetilde{z}_2)^2
+ 4 z_1 \widetilde{z}_1 z_2 \widetilde{z_2}}\,$. Parametrize $S$
away from $w_0 = 0$ by
$$
(\zeta,\eta) \mapsto [1,w_1,w_2,w_3] = [1, \zeta\eta, -\eta, \zeta] \,.
$$
The incidence relations (\ref{I}) read
\begin{eqnarray*}
z_1 - \zeta\eta \widetilde{z}_2 & = & \eta  \\
z_2 + \zeta\eta \widetilde{z}_1 & = & \zeta \;.
\end{eqnarray*}
Solving,
\be
\zeta = \frac{1+z_1 \widetilde{z}_1 + z_2 \widetilde{z}_2 \pm s}{2\widetilde{z}_2} \;.
\label{Q1} \ee
This has branch points when the square root $s$ is zero; on
$\RR^4_0$ this occurs on $C: z_1=0, |z_1|^2 + |z_2|^2 = 1$, a
circle in $\RR^3_0$. On any open set $A \subset \RR^3_0 \setminus
C$, on fixing the $\pm$ sign in (\ref{Q1}), we obtain a smooth
solution. Next note that
\begin{eqnarray*}
\theta_a & = & 2a_0 dw_1 + w_1 dw_2 - w_2 dw_1 - dw_3 \\
  & = & 2a_0(\eta d\zeta + \zeta d\eta) - \zeta\eta d\eta +
  \eta (\eta d\zeta + \zeta d\eta) - d\zeta \\
  & = & (\eta^2 + 2a_0 \eta - 1)d\zeta + 2a_0 \zeta d\eta \;.
\end{eqnarray*}  
So (\ref{E}) reads
$$
(\eta^2 + 2a_0\eta - 1) \frac{\pa\widetilde{\zeta}}{\pa\eta}
      -2a_0 \zeta \frac{\pa\widetilde{\zeta}}{\pa\zeta} = 0
$$
which has a solution
$$
\widetilde{\zeta} = \zeta \left\{ \frac{\eta + a_0 +
\sqrt{a_0^{{}2} +1}} {\eta + a_0 - \sqrt{a_0^{{}2} + 1}} \right\}
^{-a_0/\sqrt{a_0^{{}2} + 1}} \;.
$$
When $a_0 = 0$ this gives $\widetilde{\zeta} = \zeta$ so that, for
any smooth branch $\zeta:A \ra \CC$ of (\ref{Q1}) the level
curves of $\phi = \zeta|A$ give a conformal foliation. 
Explicitly,
$$
\phi(x_1,x_2,x_3) = \frac{1+x_1^2+x_2^2+x_3^2 \pm s}{2(x_2 + \ii
x_3)}
\;.
$$
with $s = \sqrt{1 -(x_1^2+x_2^2+x_3^2)^2 + 4 x_1^2 (x_2^2 +x_3^2)}
\,,$ giving a rotationally symmetric foliation.
\end{example}

Sometimes we can solve (\ref{E}) explicity even when the twistor
surface has higher degree than two; here is a cubic example:
\begin{example}
Let the twistor surface $S$ be
$$
S = \{ [w_0,w_1,w_2,w_3] \in \CP{3} : 
w_1 w_2^{{}2} + \ii w_0^{{}2} w_3 = 0 \} \;.
$$
We can parametrize this away from $w_0 = 0$ by
$$
(\zeta, \eta) \mapsto [1, z_1, z_2, z_3] = [\zeta,\ii\zeta,
\ii\eta, \ii\zeta\eta^2] \;.
$$
The incidence relations (\ref{I}) read
\be
z_1 - \ii\zeta \widetilde{z}_2 = \ii\eta \quad
z_2 + \ii\zeta \widetilde{z}_1 = \zeta \eta^2
\label{cubic-I}\,. \ee
Then
\begin{eqnarray*}
\theta_a & = & 2a_0 dw_1 + w_1 dw_2 - w_2 dw_1 - dw_3 \\
& = & (2\ii a_0 - \eta^2 + \eta )d\zeta + (-2\zeta\eta - \zeta)d\eta
\end{eqnarray*}
so (\ref{E}) reads
$$
(2\ii a_0 - \eta^2 + \eta )\frac{\pa\widetilde{\zeta}}{\pa\eta}
 + (2\zeta\eta + \zeta)\frac{\pa\widetilde{\zeta}}{\pa\zeta} = 0 \;.
$$
This can be solved by the product method to yield
$$
\widetilde{\zeta} = - \zeta \big(2\ii a_0 - \eta^2 + \eta \big)
\left( \frac{\eta-\alpha_1}{\eta-\alpha_2}
\right)^{2/\sqrt{1+8\ii a_0}}
$$
where $\alpha_1 = 1 + \sqrt{1+8\ii a_0} /2$ and
$\alpha_2 = 1 - \sqrt{1+8\ii a_0} /2$.
When $a_0 = 0$, \ $\widetilde{\zeta} = \zeta (\eta-1)^3 / \eta$
where $(\zeta, \eta)$ is a solution to (\ref{cubic-I}).  For any
$a \in \CC^4$, $\phi_a = \widetilde{\zeta}|\RR^3_a$ gives a
conformal foliation by curves.
\end{example} 


 
\leavevmode
\begin{figure}[htb]
\epsfxsize=6.5in
\centerline{\epsfbox{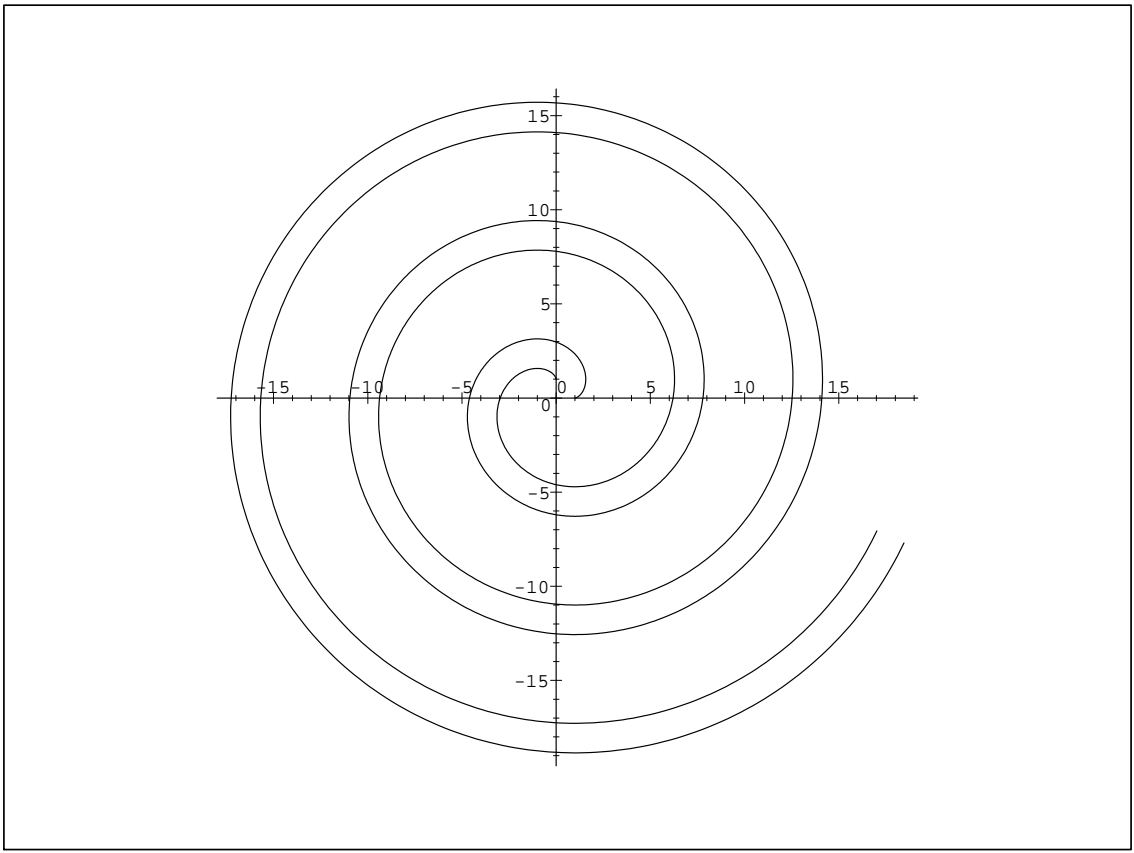}}
\caption{Two leaves of the Riemannian foliation of Example 3.10
with $t=1$.}
\end{figure}

\end{document}